\documentclass{article}

\usepackage[14pt]{extsizes}

\usepackage[utf8]{inputenc}
\usepackage[left=4cm, right = 3cm]{geometry}
\usepackage[english,russian]{babel}
\usepackage{math}
\usepackage{xcolor}
\usepackage{amssymb}
\usepackage{amsmath}
\usepackage{graphics}
\usepackage{tikz}
\usepackage{pgfplots}
\usepackage{comment}
\usepackage[ruled, linesnumbered, noend]{algorithm2e}
\title{Внешние биллиарды вне правильного десятиугольника: периодичность почти всех орбит и существование апериодической орбиты}

\author{Ф. Д. Рухович}
\begin{document}

\maketitle

%Раздел: теория динамических систем

\begin{abstract}
Доказано существование апериодической орбиты для внешнего биллиарда вне правильного десятиугольника, а также, что почти все траектории такого внешнего биллиарда являются периодическими; явно выписаны все возможные периоды.

 Работа поддержана грантом РНФ № 17-11-01337.
\end{abstract}

\begin{section}{Введение}

Пусть $\gamma$ - выпуклая фигура на плоскости $\mathbb{R}^2$, а $p$ — точка вне ее. Проведем правую относительно $p$ касательную к $\gamma$; определим $Tp \equiv T(p)$ как точку, симметричную $p$ относительно точки касания.

%\begin{comment}
\begin{Def} \label{base}
Отображение $T$ называется внешним биллиардом; фигура $\gamma$ называется столом внешнего биллиарда.
\end{Def}

Обратным к такому преобразованию является <<левый>> внешний биллиард; будем обозначать его как $T^{-1}$.

\begin{Def}
Точку $p \in \outtable$ назовем периодической, если существует такое натуральное $n$, что $T^np = p$; минимальное такое $n$ назовем периодом точки $p$ и обозначим как $per(p)$. 
\end{Def}
\begin{Def}
Точку $p \in \outtable$ назовем апериодической, если она --- не периодическая, а ее траектория бесконечна в две стороны. 
\end{Def}
\begin{Def}
Точку $p \in \outtable$ назовем граничной, если $T^np$ не определено для некоторого $n \in \bbz$.
\end{Def}

В данной статье будем полагать, что $\gamma$ --- выпуклый многоугольник.

Внешние биллиарды были введены Бернардом Нойманном в 1950-х годах и стали популярны в 1970-х благодаря Ю.Мозеру \cite{Moser78}. Внешние биллиарды исследовались рядом авторов\ (см. например, \cite{Tab93}, \cite{Schwartz09}, \cite{DF09}, \cite{SV87}, \cite{Kolodziej89}, \cite{GS91}, а также монографию \cite{Tab05}). Так, Р.Шварц \cite{Schwartz09} показал, что траектория начальной точки может быть неограниченной, тем самым разрешив вопрос Мозера\ -\ Нойманна, поставленный в \cite{Moser78}.

В центре нашего внимания находятся следующие открытые в общем случае \underline{проблемы периодичности}:

\begin{enumerate}
\item Существует ли апериодическая точка для внешнего биллиарда вне правильного $n$-угольника?
\item Верно ли, что периодические точки образуют вне стола множество полной меры для внешнего биллиарда вне правильного $n$-угольника?
\end{enumerate}

Cлучаи $n=3,4,6$ являются решеточными и тривиальными; в этих случаях, апериодической точки нет, а периодические точки, как следствие, образуют множество полной меры. По мнению Р.Шварца, следующими по сложности исследования являются случаи $n = 5,10,8,12$, ибо, по-видимому, только в этих случаях имеют место самоподобные структуры. С.Л.Табачников в \cite{Tab93} в деталях исследовал случай $n = 5$ - для них апериодические точки существуют, но их мера равна нулю. В дальнейшем правильный пятиугольник и связанная с ним символическая динамика подробно исследовались в работе N.Bedaride и J.Cassaigne \cite{BC11} (см. также их монографию \cite{BC11a}).

В статье \cite{Rukhovich18} в деталях рассматривается внешний биллиард для случая $n = 8$, а в статье \cite{Rukhovich18-2} --- для случая $n = 12$. В данной же статье исследуется случай $n = 10$. Этот случай похож на случай $n = 5$; связь между ними была доказана в \cite{BC11}. Однако возникающие в случае $n = 10$ периодические структуры, хоть и похожи на периодические структуры в случае $n = 5$, все же имеют свои особенности. Вследствие этого, детальное описание случая $n = 10$ видится необходимым.

Основным результатом данной работы являются следующие теоремы.

\begin{Th} \label{MainTreorem}
Для внешнего биллиарда вне правильного десятиугольника существует апериодическая точка.
\end{Th}

\begin{Th} \label{MainTreorem2}
В случае внешнего биллиарда вне правильного десятиугольника, периодические точки образуют вне стола множество полной меры.
\end{Th}

\begin{Th} \label{MainTheorem3}
Пусть $B_2 = \{\frac{5}{7}(6^{l+2}-(-1)^l), \frac{5}{7}(9*6^{l+1} + 2*(-1)^l), 20*6^l,
          30, 90*6^l, 10, 5, \frac{20}{7}((78+120k)*6^l - (k+1)*(-1)^l),
          \frac{5}{7}((276+240k)*6^l - (2k+3)*(-1)^l), \frac{5}{7}((234+180k)*6^l + (2k+4)*(-1)^l),
          \frac{5}{7}((34 + 40k)*6^l + (2k+1)*(-1)^l), \frac{10}{7}((20 + 40k)*6^l + (2k+2)*(-1)^l),
          40k+70, \frac{5}{7}((306 + 180k)*6^l + (2k+2)*(-1)^l), 40k+50, 60k + 40, 30k + 35,
          20k+30, 20k + 20, 10k + 15, \frac{10}{7}(6^{l+2}-(-1)^l),  \frac{10}{7}((276+240k)*6^l - (2k+3)*(-1)^l), \frac{10}{7}((34 + 40k)*6^l + (2k+1)*(-1)^l), 60k+70, 20k+30| k, l \in \bbz_{\geq 0}\}$.
          
Тогда $B_2$ есть множество всевозможных периодов периодических точек для внешнего биллиарда вне правильного десятиугольника.
% {\bf !!!!! ОПИСАТЬ ПЕРИОДЫ ДЛЯ ДЕСЯТИУГОЛЬНИКА !!!!!}
\end{Th}
 
\end{section}

\begin{section}{Внешний биллиард вне многоугольников: базовые определения и замечания} \label{basicsOfPolygonalOuterBilliard}

Во всех числовых индексах данной статьи числа будем подразумевать целыми, если не сказано обратное.

Введем основные определения, связанные с внешними биллиардами вне выпуклых многоугольников (здесь мы ссылаемся на \cite{BC11}). Пусть стол $\gamma \subset \mathbb{R}^2$ есть произвольный выпуклый $n$-угольник, $n \geq 3$. Занумеруем его вершины как $A_0A_1\ldots A_{n-1}$ против часовой стрелки. Проведем лучи $A_1A_0$, $A_2A_1$, \ldots, $A_0A_{n-1}$; они делят $\mathbb{R}^2 \backslash \gamma $ на $n$ углов, вершины которых суть вершины $\gamma$; пусть $V_i$, $0 \leq i < n$, есть один из этих углов с вершиной $A_i$. На рис. \ref{pic:baseAngleRaw1} изображен пример этих обозначений в случае $n = 5$.

\begin{figure}[h!]
\begin{center}
\includegraphics[width=90mm]{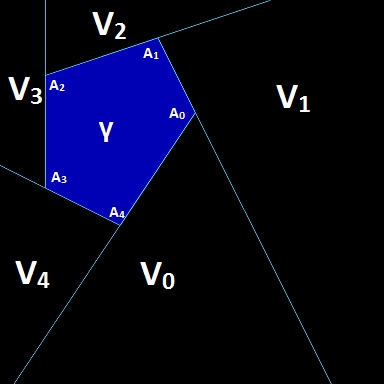}
\caption{Определение стола $\gamma$, вершин и углов $V_i$}.
\label{pic:baseAngleRaw1}
\end{center}
\end{figure}

%{\bf !!!! ДОБАВИТЬ РИСУНОК ПРО $V_i$ и Ко !!!!}

Из определения внешнего биллиарда напрямую следует

\begin{Lm} \label{basicAngleV}
Пусть $p \in \outtable$, а $i \in [0, n)$. Тогда преобразование $T$ для точки $p$ определено и является центральной симметрией относительно вершины $A_i$, если и только если $p \in \intt (V_i)$.
\end{Lm}

Поймем, как выглядит множество граничных точек.

\begin{Lm} \label{boundaryPointsAreOpenRaysAndSegments}
Множество граничных точек есть объединение счетного числа открытых отрезков и лучей. 
\end{Lm}

\begin{Proof}
Если удалить из прямой, являющейся продолжением стороны $\gamma$, отрезок, являющийся стороной $\gamma$, то прямая распадется на два открытых луча. Проделаем такую операцию для каждой из $n$ сторон; пусть $L$ есть множество точек, лежащих на таким образом полученных $2n$ открытых лучах. Очевидно, что множество граничных точек есть $\{T^k(p) | p \in L, k \in \bbz \}$, т.е. результат применения неограниченного числа преобразований $T$ и $T^{-1}$ к лучам, образующим $L$. Однако заметим, что если $T$ (или $T^{-1}$) определено для хотя бы одной точки некоторого отрытого луча или отрезка, то $T$ (или $T^{-1}$) разделит этот луч или отрезок лучами $A_{(i+1)\ mod\ n}A_i$, $0 \leq i < n$, на конечное число открытых отрезков и/или лучей, для каждого из которых $T$ (или $T^{-1}$) есть центральная симметрия. Тогда для каждого $k \in \bbz$, $T^k(L)$ есть объединение конечного числа отрезков и лучей, а множество граничных точек есть объединение счетного числа отрезков и лучей как объединение конечных множеств отрезков и лучей, QED.
\end{Proof}

Так прямая, луч и отрезок являются в $\bbrr$ множествами меры нуль, то из леммы \ref{boundaryPointsAreOpenRaysAndSegments} 
напрямую следует известный факт, важный ввиду второй проблемы периодичности.

\begin{Lm} \label{zeroMeasureOfBoundaryPoint}
Множество граничных точек для внешнего биллиарда вне $\gamma$ является множеством меры нуль. 
\end{Lm}

Введем <<классическое>> для внешнего биллиарда вне правильных многоугольников {\it кодирование}.

\begin{Def}
Пусть $p \in \outtable$ --- периодическая или апериодическая точка. Тогда кодом $\rho(p) \equiv \rho_{\gamma}(p)$ является последовательность ($\ldots u_{-2}u_{-1}u_0u_1u_2\ldots $), т.ч. $\forall i \in \bbz: T^i(p) \in \intt(V_{u_i})$.
\end{Def}

\begin{Def}
Пусть $p \in \outtable$ --- граничная точка,последовательное применение внешнего биллиарда $T$ может быть выполнено ровно $m$ раз, а преобразования $T^{-1}$ - $l$ раз, $0 \leq l, m \leq +\infty$. Тогда код $\rho(p) = \rho_{\gamma}(p)$ есть последовательность $(u_i)_{i \in [-l, m)}$, где $u_i \in [0, n)$ и $u_i = k$, если и только если $T^i(p) \in \intt(V_k)$.
\end{Def}

Также будем обозначать элемент $u_i$ кода как $\rho(p)[i]$, а подпоследовательность $u_lu_{l+1}\ldots u_r$, $-\infty < l \leq r < +\infty$ --- как $\rho(p)[l, r]$.

\begin{Def}
Компонентой назовем максимальное по включению множество точек с одинаковым кодом $\rho$; компоненту, в которой содержится точка $p$, обозначим за $comp(p)$.
\end{Def}

Непосредственно из определений следует

\begin{Lm} \label{componentCalculation}
Пусть точка $p$ обладает кодом $\rho(p)$ бесконечной длины, равным $\ldots u_{-2}u_{-1}u_0u_1u_2 \ldots$. Тогда $comp(p) = \underset{i \in \bbz}{\cap} U_i$, где $U_i$ есть:

\begin{itemize}
    \item $\intt(V_{u_0})$, при $i = 0$;
    \item $T^{-i}(\intt(V_{u_i}) \cap T^{i}(U_{i-1}))$, при $i > 0$;
    \item $T^{-i}(\intt(V_{u_i}) \cap T^{i}(U_{i+1})$, при $i < 0$.
\end{itemize} 
\end{Lm}

Отметим, что множество $U_i, i \in \bbz$ представляет собой пересечение конечного числа полуплоскостей, площадь которого отлична от нуля, но может быть бесконечной. Если пересечение ограничено, то это пересечение есть выпуклый многоугольник. Неограниченное пересечение полуплоскостей будем называть бесконечным многоугольником; если же пересечение ограничено, назовем его конечным многоугольником. 

Поймем, как выглядит компонента периодической точки.

\begin{Lm} \label{periodic1}
Пусть $p$ - периодическая точка c периодом $m$. Пусть $q$ - некоторая точка, т.ч. определено $T^{2m}(q)$, причем $\rho(p)[0, 2m-1] = \rho(q)[0, 2m-1]$. Тогда:

\begin{enumerate}
\item $q$ - периодическая точка с периодом не более чем $2m$;
\item $q \in comp(p)$, т.е. $\rho(p) = \rho(q)$.
\end{enumerate}
\end{Lm}

\begin{Proof}
Пусть $\vec v = \vec {pq}$. Тогда в силу одинаковости соответствующих частей кода и по свойствам центральной симметрии $T^k(q) = T^k(p) + (-1)^k \vec v$, $k = 0, 1, \ldots, 2m$; в частности, $T^{2m}(q) = T^{2m}(p) + (-1)^{2m} \vec v = p + v = q$. Таким образом, $q$ периодична, а последовательность $T^k(q)$ имеет (возможно, не минимальный) период $2m$. Следовательно, $\rho(q)$ бесконечен в обе стороны, имеет период $2m$ и, как очевидное следствие, совпадает с $\rho(p)$.
\end{Proof}

Прямым следствием предыдущей леммы является
\begin{Lm} \label{allInCompPeriodic}
Пусть $p$ - периодическая точка c периодом $m$. Тогда:
\begin{itemize}
    \item $comp(p) = U_{2m}$ (определение множества $U$ см. в условии леммы \ref{componentCalculation});
    \item все точки $comp(p)$ периодические, причем каждая из них обладает (возможно, не минимальным) периодом $2m$;
    \item $comp(p)$ есть открытый, конечный или бесконечный, выпуклый многоугольник, стороны которого параллельны сторонам $\gamma$.
\end{itemize}
\end{Lm}

Последнее утверждение леммы \ref{allInCompPeriodic} очевидно следует из структуры множеств $U_{i}$. Понять же устройство периодической компоненты, т.е. компоненты, содержащей периодическую точку, в бОльших деталях нам поможет следующая лемма.

\begin{Lm} \label{diffCodeDividing}
Пусть $p$, $q$ - две точки вне стола $\gamma$, и пусть для некоторого $l \in \bbz$ $\rho(p)_l$ и $\rho(q)_l$ определены, причем $\rho(p)_l \neq \rho(q)_l$. Тогда на отрезке $pq$ существует граничная точка.
\end{Lm}

\begin{Proof}
Докажем лемму для случая $l \geq 0$ (случай $l < 0$ может быть рассмотрен аналогично). Без ограничения общности будем считать, что $l$ минимально среди всех подходящих $l$. Тогда для любой точки $x$ отрезка $pq$ $T^l(x)$ определено и является последовательностью центральных симметрий относительно точек $A_{u_0}$, $A_{u_1}$, ..., $A_{u_{l-1}}$, где $\ldots u_0u_1u_2\ldots = \rho(p)$. По свойствам центральных симметрий, $T^l(pq)$ есть отрезок $T^l(p)T^l(q)$. Так как $\rho(p)_l \neq \rho(q)_l$, то $T^l(p)$ и $T^l(q)$ лежат в разных углах $V_{i_p}$ и $V_{i_q}$; следовательно, отрезок $T^l(p)T^l(q)$ пересекает границу угла $V_{i_p}$ (равно как и $V_{i_q}$) в некоторой точке $y$, а точка $x = T^{-l}(y)$ лежит на отрезке $pq$ и является граничной точкой, QED.
\end{Proof}

Такая лемма дает возможность понять, из чего состоит граница компоненты периодической точки.

\begin{Lm} \label{allOnBoundFinite}
Пусть $p$ - периодическая точка. Тогда $\partial comp(p)$ состоит лишь из граничных точек (точки самого стола $\gamma$ будем также считать граничными).
\end{Lm}

\begin{Proof}
Рассмотрим точку $q$, лежащую на границе $comp(p)$. Предположим, что $q$ не является граничной; тогда $q$ обладает бесконечным в обе стороны кодом $\rho(q)$. Так как $q \notin comp(p)$, то $\rho(p) \neq \rho(q)$; следовательно, по лемме \ref{diffCodeDividing} на отрезке $pq$ существует граничная точка с конечной траекторией. Полученное противоречие с леммой \ref{allInCompPeriodic} завершает доказательство.
\end{Proof}

\begin{Lm} \label{periodicComponentIsConvexPolygon}
Пусть $p$ - периодическая точка. Тогда $comp(p)$ есть конечный многоугольник.
\end{Lm}

\begin{figure}[h!]
\begin{center}
\includegraphics[width=150mm]{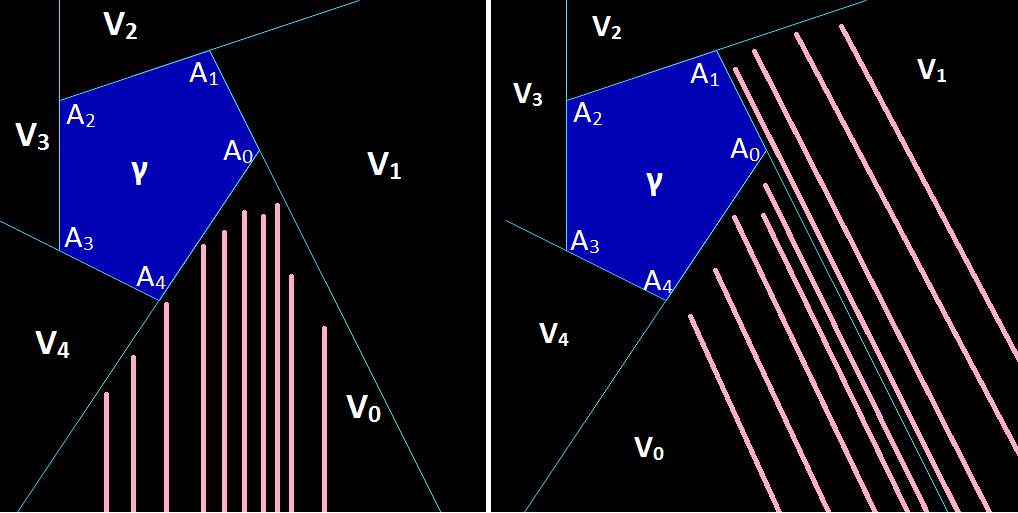}
\caption{Сонаправленные лучи, лежащие целиком в одном из секторов}.
\label{pic:baseAngleParallelRays}
\end{center}
\end{figure}

\begin{Proof}
Пусть это не так. Тогда по лемме \ref{allInCompPeriodic} $comp(p)$ есть бесконечный выпуклый многоугольник; следовательно, существует луч $r$ с началом в некоторой точке $c$ и направляющим вектором $\vec v$, целиком лежащий внутри $comp(p)$.

Из свойств центральной симметрии следует, что $T^{2l}(r)$, $l \in \bbz_{\geq 0}$ есть луч, сонаправленный с лучом $r$. Каждый из  лучей $T^{2l}(r)$ должен лежать целиком внутри одного из углов $V_i$. Заметим (см. рис. \ref{pic:baseAngleParallelRays}), что если рассмотреть всевозможные лучи с заданным направлением $\vec v$, каждый из которых лежит целиком строго внутри какого-то из углов $V_i$, то эти лучи либо лежат в одном из углов, либо в двух соседних углах, если луч сонаправлен с лучом - стороной одного из углов.

Докажем, что последовательность лучей $T^{2l}(r)$ лежит строго в одном из углов. Пусть это не так, и лучи 
$T^{2l}(r)$ лежат в двух соседних углах; пусть это, без ограничения общности, $V_0$ и $V_1$ (нумерация углов, напомним, идет против часовой стрелки). Так как $T^{2l}(r)$ есть периодическая последовательность (доказательство можно провести аналогично лемме \ref{periodic1}), то $\exists l_0 \in \mathbb{N}: T^{2l_0}(r) \subset \intt(V_0) \wedge T^{2l_0 + 2}(r) \subset \intt(V_1)$; следовательно, существует луч $r' \subset \intt(V_0)$, сонаправленный с $r$, т.ч. $T^2(r') \subset \intt(V_1)$. По свойствам центральных симметрий, для $r'$ $T^2$ есть параллельный перенос на вектор $2 \vec u$, где $\vec u$ - вектор, начинающийся в вершине $A_0$ угла $V_0$ и заканчивающийся в некоторой другой вершине стола $\gamma$. Т.к. $\gamma$ лежит в левой полуплоскости относительно ориентированной прямой $A_0A_1$, то и вектор $\vec u$ направлен <<невправо>> относительно прямой $A_0A_1$. Как следствие, лучи $r'$ и $T^2(r')$ лежат в левой относительно $A_0A_1$ полуплоскости, а угол $\intt(V_1)$ лежит в правой относительно $A_0A_1$ полуплоскости - противоречие.

Итак, все лучи последовательности $T^{2l}(r)$, $l \in \bbz$, лежат в некотором угле $V_i$, а, аналогично, лучи последовательности $T^{2l+1}(r)$, $l \in \bbz$ - в некотором угле $V_j$. Тогда код всех вершин луча есть $\ldots ijijijijijijijijij \ldots $; следовательно, по все тем же свойствам центральных симметрий, $T^{2l}(c) = c + 2l\vec u$, $\vec u = \overrightarrow{A_iA_j}$. Следовательно, $\forall l \in \bbz: T^{2l}(c) \neq c$ - противоречие с периодичностью точки $p$ относительно $T$ и, как следствие, $T^2$. Таким образом, $comp(p)$ есть ограниченное множество, QED.
\end{Proof}

Итак, компонента периодической точки $p$ есть открытый выпуклый не-более-чем-$2n$-угольник, стороны которого параллельны сторонам стола $\gamma$. Более того, из лемм \ref{allInCompPeriodic}, \ref{allOnBoundFinite} вытекает альтернативное определение понятия <<компонента периодической точки>>:

\begin{Lm} \label{periodicCompIsConnectionComp}
Пусть $p \in \outtable$ - периодическая точка. Тогда компонента периодической точки $comp(p)$ есть максимальное по включению связное множество периодических точек, содержащее $p$.
\end{Lm}

Но как же устроены периоды точек $comp(p)$?

\begin{Lm} \label{commonPeriod}
Пусть $p$ - периодическая точка. Тогда $\forall n \in \mathbb{N}: T^n(comp(p)) = comp(p) \vee T^n(comp(p)) \cap comp(p) = \varnothing$.  
\end{Lm}

\begin{Proof}
$T^n(comp(p))$ есть открытый выпуклый многоугольник, равный $comp(p)$; следовательно, ни одно их множеств $T^n(comp(p))$ и $comp(p)$ не может быть строгим подмножеством другого множества. В этом случае, либо вышеописанные множества не пересекаются, либо совпадают, либо пересекаются таким образом, что существует точка $q$, которая лежит внутри $T^n(comp(p))$ и на границе $comp(p)$. Последний случай вызывает противоречие, ибо по леммам \ref{allInCompPeriodic}, \ref{allOnBoundFinite} точка $q$ должна быть и периодической, и граничной. Следовательно, лемма доказана.    
\end{Proof}

Предыдущая лемма дает возможность ввести понятие периода компоненты.

\begin{Def}
Пусть $p$ - периодическая точка. Тогда периодом компоненты $comp(p)$, или $per(comp(p))$, назовем минимальное такое натуральное $k$ такое, что $T^k(p) = p$.
\end{Def}

Как же $per(comp(p))$ связано с периодами точек $comp(p)$?

Следующая лемма есть прямое следствие леммы \ref{commonPeriod} и того факта, что $T^n$ для компоненты периодической точки $p$ есть параллельный перенос в случае четного $n$ и центральная симметрия в случае нечетного $n$.

\begin{Lm} \label{fixedPeriods}
Пусть $p$ - периодическая точка, и пусть $k = per(comp(p))$. Тогда:
\begin{itemize}
\item если $k$ четно, то все точки $comp(p)$ имеют период $k$;
\item если $k$ нечетно, то: а) $comp(p)$ есть центрально-симметричный многоугольник с центром в некоторой точке $c$; б) период точки $c$ есть $k$, а всех остальных точек - $2k$.
\end{itemize}
\end{Lm}

Пользуясь данной леммой, будем также называть компоненту некоторой периодической точки {\it периодической компонентой}.

Согласно лемме \ref{fixedPeriods}, мы можем найти множество всевозможных периодов точек для внешнего биллиарда вне $\gamma$, если нам известно множество всевозможных периодов периодических компонент. Точное соотношение между множествами задает следующая лемма, являющаяся прямым следствием леммы \ref{fixedPeriods}.

\begin{Lm} \label{pointsPeriodsFromComponentsPeriods}
Пусть $B_c$ и $B_p$ --- множества всевозможных периодов компонент и точек соответственно для внешнего биллиард вне многоугольника $\gamma$. Тогда $B_p = B_c \cup \{2 * l\ |\ l \in B_c,\ l \text{ нечетно}\}$.
\end{Lm}

Завершит наше исследование на тему периодических компонент следующая лемма.

\begin{Lm} \label{invariancyOfPeriodicComponent}
Пусть $p$ - периодическая точка. Тогда $T(comp(p)) = comp(T(p))$.
\end{Lm}
\begin{Proof}
Очевидно, что $T(comp(p)) \subset comp(T(p))$ и $T^{-1}(comp(T(p))) \subset comp(T^{-1}(T(p))) = comp(p)$, откуда следует, что $comp(T(p)) \subset T(comp(p))$; следовательно, $T(comp(p)) = comp(T(p))$, QED.
\end{Proof}

\end{section}

\begin{section}{Внешний биллиард вне правильного десятиугольника}

Перейдем к исследованию правильного десятиугольника. С этого момента будем считать, что стол внешнего биллиарда есть правильный десятиугольник $\gamma = A_0A_1 \ldots A_9$, вершины которого перенумерованы против часовой стрелки. 
 Будем следовать плану, намеченному в \cite{Rukhovich18, Rukhovich18-2}.

\begin{subsection}{Ограничение преобразования}

Введем ограничение преобразования $T$, похожее на ограничение, выполненное в работах \cite{BC11, Tab93}, но отличающееся в итоговой реализации.

\begin{Def}
Пусть $R$ есть поворот на угол $\frac{\pi}{5}$ по часовой стрелке вокруг центра $\gamma$.
\end{Def}

Заметим, что $T$ инвариантно относительно $R$, т.е. $\forall p \in \outtable$, $T(p) \text{ определено: } T(R(p)) = R(T(p))$. Отождествим точки, переходящие друг в друга с помощью $R$.

Пусть $V'$ есть угол с вершиной в $A_1$, центрально-симметричный углу $V_1$ относительно $A_1$.

\begin{Def}
Пусть $p \in \outtable$. Тогда $k_p$ есть такое минимальное число $k \in \bbz_{\geq 0}$, что $R^k(p) \in V'$, a $R'(p)$ есть $R^{k_p}(p)$.
\end{Def}

Другими словами, $R'(p)$ есть представитель класса эквивалентности $p$ в $V'$ (если только $p$ не лежит на продолжении одной из сторон).

\begin{Def}
%Пусть $V'$ есть угол с вершиной в т. $A_1$, образованный лучами $A_0A_1$ и $A_1A_2$. Тогда отождествление индуцирует преобразование $T': V' \rightarrow V'$

Преобразование $T': V' \rightarrow V'$ есть преобразование, индуцированное отождествлением точек относительно $R$, устроенное следующие образом. Пусть $p \in \intt(V')$; тогда:

\begin{itemize}
    \item $T'(p)$ определено, если и только если $T(p)$ определено;
    \item если $T'(p)$ определено, то $T'(p) = R'(T(p))$.
\end{itemize}

Для удобства будем считать, что на $\partial V'$ $T'$ не определено.

\end{Def}

Любую точку $p \in V'$ будем рассматривать как граничную (периодическую, апериодическую) относительно $T'$ ровно в том же смысле, в котором ранее мы рассматривали точки $p \in \outtable$.
Из структуры отождествления очевидны следующая лемма.

\begin{Lm} \label{TandInducedT}
  Точка $p \in V'$ является граничной (периодической, апериодической) относительно $T'$, если и только если $p \in V'$ является граничной (периодической, апериодической) относительно $T$. Более того, в этом случае граничными (периодическими, апериодическими) относительно $T$ будут являться точки $R^k(p)$, $k = 0, 1, \ldots, 9$.
\end{Lm}
 
 Прямым следствием леммы $\ref{TandInducedT}$ является следующая лемма, которая сводит решение проблем периодичности для всей плоскости относительно преобразования $T$ к решению тех же проблем, но в $\intt(V')$ и относительно преобразования $T'$.

\begin{Lm} \label{reductionToT'}
\begin{enumerate}
    \item Апериодическая относительно преобразования $T$ точка существует, если и только если существует апериодическая точка в $\intt(V')$ относительно преобразования $T'$;
    
     \item Периодические относительно $T$ точки образуют вне $\gamma$ множество полной меры, если и только если периодические относительно $T'$ точки образуют в $\intt(V')$ множество полной меры. 
\end{enumerate}
    
\end{Lm}    

Чтобы свести проблему нахождения периодов к преобразованию $T'$, 
введем индуцированный ограничением код относительно $T'$.

\begin{Def}
Пусть $p \in V'$ --- периодическая или апериодическая точка. Тогда индуцированным кодом $\rho'(p) \equiv \rho'_{\gamma}(p)$ является последовательность ($\ldots v_{-2}v_{-1}v_0v_1v_2\ldots $), т.ч. $\forall i \in \bbz: T'^i(p) \in \intt(V_{v_i + 1})$.
\end{Def}

Заметим, что коды $\rho$ и $\rho'$ связаны между собой следующим образом.

\begin{Lm} \label{codesConnection}
Пусть $p \in V'$ --- периодическая или апериодическая точка. Тогда $\forall i \in \bbz$: $\rho'(p)[i] = (\rho(p)[i] - \rho(p)[i-1]) \mod 10$.
\end{Lm}

\begin{Proof}
Преобразования $T$, $T'$ и $R'$ устроены таким образом, что если $T(q_1) = q_2$, то $T'(R'(q_1)) = R'(q_2)$, и $R'(R(q_1)) = R'(q_1)$. Отсюда заметим, что если $q$ - произвольная неграничная точка, то при замене $q$ на $R(q)$ все значения кода $\rho(q)$ уменьшаются на один по модулю 10, а значения кода $\rho'(R'(q))$ не изменяются. Зафиксируем произвольное $i \in \bbz$. Пусть $k = k_{T^i(p)}$, и $p' = R^k(p)$. Тогда $\rho'(R'(p'))[i] = \rho'(p)[i]$, и $(\rho(p')[i] - \rho(p')[i-1]) \mod 10\ =\ (\rho(p)[i] - \rho(p)[i-1]) \mod 10$. С другой стороны, $T^i(p') \in \intt(V')$, а $T^{-1}(\intt(V')) = \intt(V_1)$. Следовательно, по определению $\rho(p')[i-1] = 1$, $\rho(p')[i] = 1 + \rho'(R'(p'))[i] = \rho(p')[i-1] + \rho'(p)[i]$, откуда получаем $(\rho(p)[i] - \rho(p)[i-1]) \mod 10 = \rho'(p')[i]$. Факт произвольности выбора $i$ завершает доказательство леммы. 
\end{Proof}

Очевидно, что внутри $V'$ периодические компоненты относительно $\rho$ и $\rho'$ одни и те же (ибо для любой неграничной точки $p \in V'$ код $\rho(p)$ можно восстановить из $\rho'(p)$ и наоборот), что дает возможность ввести период $per'(comp)$ как период периодической компоненты $comp \subset V'$ относительно $T'$. Отметим, что все 
элементы $\rho'(p)$ являются целыми числами от 1 до 5.

\begin{Lm} \label{fromInducedPeriodToPeriod}
Пусть $p \in V'$ - периодическая точка, причем $per'(p) = m$. Тогда $per(p) = m * \frac{10}{\rgcd(s, 10)}$, где $s = \sum\limits_{i=1}^{m} \rho'(p)[i]$.
\end{Lm}

\begin{Proof}
Заметим, что для любого $k \in \bbz$ верна импликация: $(T^k(p) = p) \implies (T'^k(p) = p)$; следовательно, $per(p)$ кратен $m$; пусть $z = \frac{per(p)}{m}$. По лемме \ref{codesConnection} и в силу периодичности $p$ верно: $s \mod 10 = (\rho(p)[m] - \rho(p)[0]) \mod 10 = (\rho(p)[2m] - \rho(p)[m]) \mod 10 = (\rho(p)[3m] - \rho(p)[2m]) \mod 10 = \ldots$. Так как $p = T'^m(p) = T'^{2m}(p) = \ldots$, то $\forall l \in \bbz$: $f'(T^{lm}(p)) = p$; следовательно, $\forall l \in \bbz$: $T^{lm}(p) = p \Leftrightarrow (\rho(p)[lm] - \rho(p)[0]) \mod 10 = 0$, а $z$ есть минимальное натуральное подходящее $l$. Так как $\forall l \in \bbz: (\rho(p)[lm] - \rho(p)[0]) \mod 10 = (l*(\rho(p)[m] - \rho(p)[0])) \mod 10 = (l*s) \mod 10$, то $z$ есть минимальное натуральное число такое, что $z*s$ делится на 10; из теории чисел очевидно, что $z = \frac{10}{\rgcd(s, 10)}$, QED.
\end{Proof}

Отметим, что лемма \ref{fromInducedPeriodToPeriod} останется верна, если вместо точки $p$ рассмотреть её периодическую компоненту. Таким образом, зная всевозможные {\it индуцированные коды периодов} относительно $T'$, можно описать и всевозможные периоды относительно $T$. Более того, для этого не нужно знать именно слова - достаточно знать, сколько раз в коде периодической траектории появляются единица, двойка, \ldots, пятерка. Более формально эту мысль описывает следующая лемма, являющаяся прямым следствием леммы \ref{fromInducedPeriodToPeriod}.

\begin{Lm} \label{fromIPStructureToPeriod}
Пусть $p \in (\subset) V'$ - периодическая точка (компонента), и пусть $per'(p) = m$. Пусть в последовательности $\rho'(p)[1, m]$ число $j$ встречается ровно $a_j$ раз, $j = 1, 2, \ldots, 5$. Тогда $per(p) = \frac{10*(a_1 + a_2 + \ldots + a_5)}{\rgcd(10, a_1 + 2a_2 + \ldots + 5a_5)}$.
\end{Lm}

\end{subsection}

\begin{subsection}{$T'$ как <<кусочное движение>>}.

Опишем, как выглядит преобразование $T'$. Для этого определим углы $V'_i$, $0 \leq i < 10$, как угол, симметричный $V_i$ относительно вершины $A_i$. В частности, $V' = V'_1$. Из определения очевидно, что $\forall p \in \outtable, j \in [0, 10)$: $\rho(p)[0] = j \Leftrightarrow T(p) \in \intt(V'_i)$ (в частности, $\rho(p)[0]$ определено, если и только если $T(p)$ определено), а также $\forall j \in [0, 10)$: $R(V'_j) = V'_{(j - 1) \mod 10}$. 

\begin{Lm} \label{inducedRotation}
Пусть $p \in \intt(V')$, и $T(p)$ определено. Тогда $T'(p) = R^{v_0}(T(p))$, где $v_0 = \rho'(p)[0]$.  
\end{Lm}

\begin{Proof}
По определению $\rho'$, $T$ для $p$ есть центральная симметрия относительно вершины $A_{v_0 + 1}$. Тогда $T(p) \in \intt(V'_{v_0 + 1})$. Так как $R^{v_0}(V'_{v_0 + 1}) = V'_1 = V'$ (очевидно), то $R^{v_0}(T(p)) \in \intt(V')$. Следовательно, $T'(p) = R^{v_0}(T(p))$, QED.  
\end{Proof}

Пусть $\alpha_i \subset V'$, $i = 1,2,\ldots, 5$, есть множество точек $p \in V'$, т.ч. $\rho'(p)[0] = i$; все фигуры $\alpha_i$ изображены на рис. \ref{pic:inducedT12}. Из леммы \ref{basicAngleV} очевидно, что $\alpha_i = \intt(V') \cap \intt(V_{i + 1})$, откуда следует следующая лемма.

\begin{figure}[h!]
\begin{center}
\includegraphics[width=150mm]{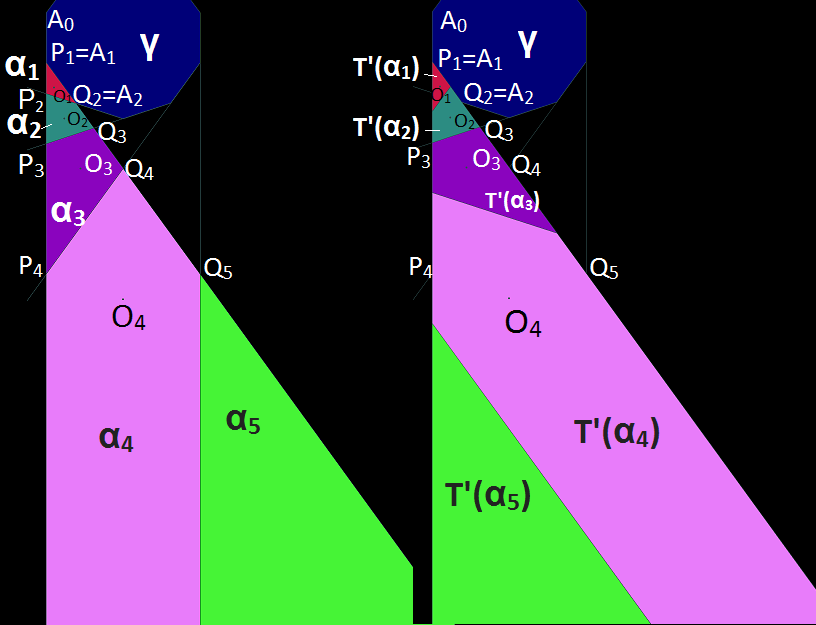}
\caption{Фигуры $\alpha_i$ и индуцированное преобразование $T'$.}
\label{pic:inducedT12}
\end{center}
\end{figure}

\begin{Lm} \label{inducedTAlphas}
\begin{enumerate}
    \item $\alpha_1$ есть открытый треугольник $P_1P_2Q_2$;
    \item $\alpha_i$, $i = 2, 3$, есть открытый четырехугольник $P_iQ_iQ_{i+1}P_{i+1}$;
    \item $\alpha_4$ есть открытый бесконечный многоугольник, ограниченный лучом $A_0A_1$, отрезками $P_4Q_4$ и $Q_4Q_5$ и лучом $A_6A_5$;
    \item $\alpha_5$ есть открытый угол между лучами $A_1A_2$ и $A_6A_5$.
\end{enumerate}
\end{Lm}

Важными для нашего описания оказываются точки $O_i$ пересечения биссектрисы угла $P_2P_1Q_2$ с биссектрисами углов $V_{i+1} = P_iA_{i+1}P_{i+1}$, $i = 1,2,3,4$.

\begin{Lm} \label{fixedPoints}
$\forall i \in \{1,2,3,4\}$: $T'(O_i) = O_i$.
\end{Lm}
\begin{Proof}
Пусть $i \in \{1,2,3,4\}$, а $O$ - центр $\gamma$. Заметим, что угол $V_{i+1}$ симметричен углу $V'$ относительно прямой $OP_{i+1}$; следовательно, $A_1O_i = A_{i+1}O_i$. С другой стороны, по лемме \ref{inducedRotation} $T'(O_i) = R^i(T(O_i))$, причем преобразование $T$ переводит с сохранением расстояния биссектрису угла $V_{i+1}$ в биссектрису угла $V'_{i+1}$, а $R^i$ переведет последнюю биссектрису с сохранением расстояния в биссектрису угла $V'_1 = V'$. Следовательно, $T'(O_i)$ будет лежать на биссектрисе угла $V'$, причем $A_1T'(O_i) = A_{i+1}O_i = A_1O_i$, откуда $O_i = T'(O_i)$, QED.
\end{Proof}

Лемма \ref{fixedPoints} позволяет нам полностью описать преобразование $T'$.

\begin{Lm} \label{inducedTAction}
\begin{enumerate}
    \item $\forall i \in \{1,2,3,4\}$: $T'(\alpha_i)$ есть поворот на угол $\frac{(5-i)\pi}{5}$ против часовой стрелки вокруг точки $O_i$;
    \item $T'(\alpha_5)$ есть параллельный перенос вдоль вектора $\overrightarrow{A_6A_1}$;
    \item $\forall i \in \{1,2,3,4,5\}~\forall p \in \partial(\alpha_i)$: $T'(p)$ не определено.
\end{enumerate}
\end{Lm}

\begin{Proof}
Пункт 3 очевидно следует из определения $T'$ и того факта, что границы $\alpha_i$ лежат на продолжениях сторон стола $\gamma$.

Докажем пп.1,2. Согласно лемме \ref{inducedRotation}, $T'(\alpha_i)$, $i = 1,2,\ldots,5$ есть $R^{i}(T(\alpha_i))$; т.к. центральная симметрия есть поворот на угол $\pi$ против часовой стрелки (будем для удобства использовать такое направление), то при $i = 1,2,3,4$, $T'(\alpha_i)$ есть поворот на угол $\pi - i*\frac{\pi}{5}$ против часовой стрелки вокруг некоторой неподвижной точки, которой по лемме \ref{fixedPoints}, и п.1 доказан. Что до п.2, то при $i = 5$ $R^i$ есть центральная симметрия относительно т. $O$, центра $\gamma$; следовательно, $T'(\alpha_5)$ есть композиция центральных симметрий относительно точек $A_6$ и $O$; из школьной геометрии известно, что это есть параллельный перенос на вектор $2{\overrightarrow{A_6O}} = {\overrightarrow{ A_6A_1}}$, QED.

\end{Proof}

Таким образом, $T'$ разбивает $V'$ на пять фигур, для каждой из которых $T'$ есть движение, т.е. $T'$ кусочно изометрично (piecewise isometric).

Введем еще одно определение.

\begin{Def}
Пусть $\alpha, \beta \subset \bbrr$ - конечные или бесконечные многоугольники. Будем говорить, что $\beta$ вписан в $\alpha$, если $\beta \subset \alpha$ и каждая из сторон $\alpha$ целиком содержит одну из сторон $\beta$.
\end{Def}

Рассмотрим периодические компоненты $comp(O_i)$, $i = 1,2,3,4$. Каждая такая компонента имеет код $\rho'(comp(O_i))$, равный $\ldots iiiiii \ldots$ и является, согласно лемме \ref{allInCompPeriodic}, открытым выпуклым многоугольником; назовем эти многоугольники $\beta_i$. Будем также обозначать отрезок с концами в точках $p, q$ как $\overline{pq}$.

\begin{Lm} \label{period1Components}
\begin{enumerate}
    \item $\beta_1$ есть правильный пятиугольник $B^1_0B^1_1B^1_2B^1_3B^1_4$, т.ч. $B^1_1 = P_2$, $\overline{B^1_0B^1_1} \subset \overline{P_1P_2}$, $\overline{B^1_1B^1_2} \subset \overline{P_2Q_2}$,
    $\overline{B^1_3B^1_4} \subset \overline{Q_2Q_1}$;
    
    \item $\beta_2$ есть правильный десятиугольник $B^2_0B^2_1 \ldots B^1_9$, т.ч. $B^2_7 = Q_2$,
    $\overline{B^2_0B^2_1} \subset \overline{P_2P_3}$, 
    $\overline{B^2_3B^2_4} \subset \overline{P_3Q_3}$, 
    $\overline{B^2_6B^2_7} \subset \overline{Q_3Q_2}$,
    $\overline{B^2_7B^2_8} \subset \overline{Q_2P_2}$;
    
    \item $\beta_3$ есть правильный пятиугольник $B^3_0B^3_1 \ldots B^3_4$, т.ч. $B^3_0 = P_3$, $B^4_3 = Q_4$, $B^3_3 = Q_5$, 
    $B^3_1 \in \overline{P_3P_4}$, 
    $B^3_2 \in \overline{P_4Q_4}$; 

    \item $\beta_4$ есть правильный десятиугольник $B^4_0B^4_1 \ldots B^4_9$, т.ч. $B^4_0 = P_4$, $B^4_6 = Q_5$, $B^4_7 \in \overline{Q_5Q_4}$,  $B^4_9 \in \overline{Q_4P_4}$.

\end{enumerate}

\end{Lm}

\begin{figure}[h!]
\begin{center}
\includegraphics[width=150mm]{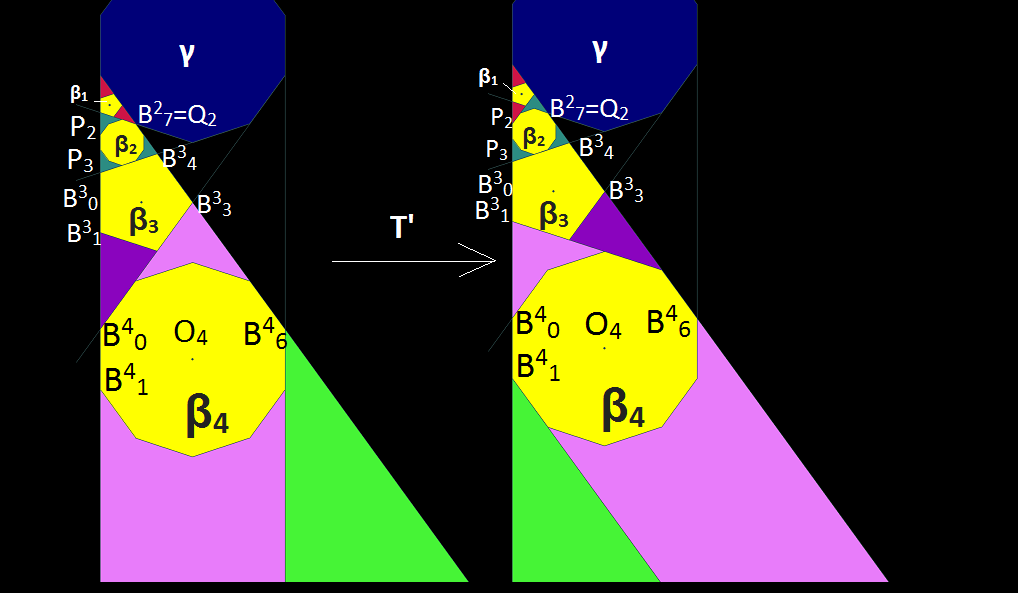}
\caption{Инвариантные фигуры $\beta_i$}
\label{pic:InducedT-InvariantComponents-12}
\end{center}
\end{figure}

Таким образом, $\beta_j$ есть правильные многоугольники, вписанные в $\alpha_j$, $j=1,2,3,4$; этим многоугольники изображены на рис. \ref{pic:InducedT-InvariantComponents-12}. Для доказательства достаточно проверить, что при применении преобразования $T'$ фигуры как открытые многоугольники переходят сами в себя, а их границы состоят из граничных точек; проверку этих фактов с помощью леммы \ref{inducedTAction} (и рис. \ref{pic:InducedT-InvariantComponents-12}) оставляем читателю.

\end{subsection}
\begin{subsection}{Самоподобие 1}

Рассмотрим угол $\alpha_5$.

\begin{Def} \label{bigFirstReturnMap}
Определим преобразование $T'': \alpha_5 \rightarrow \alpha_5$ следующим образом. Пусть $q \in \alpha_5$, и пусть $k = k_{alpha_5}(q)$ есть минимальное целое положительное число такое, что $q' = T'^k(q)$ определено, и $q' \in \alpha_5$. Тогда $T''(q) = q'$.
\end{Def}

Другими словами, $T''$ есть {\it преобразование первого возвращения} (first return map) относительно $T'$ на $\alpha_5$. 

Оказывается, $T'$ для $\intt(V')$ идентично $T''$ для $\alpha_5$! Для формального установления этого соответствия, введем преобразование $H: \intt(V') \rightarrow \alpha_5$ как параллельный перенос на вектор $A_1Q_5$.

\begin{Lm} \label{selfSimilarityBig}
Пусть $p \in \intt(V')$, $q = H(p)$. Тогда:
\begin{itemize}
    \item $T'(p)$ определено, если и только если $T''(q)$ определено;
    \item если $T'(p)$ определено, то $H(T'(p)) = T''(q) = T''(H(p))$.
\end{itemize}
\end{Lm}

Утверждение леммы иллюстрирует рис. \ref{pic:bigSelfSimilarity-12}.

\begin{figure}[h!]
\begin{center}
\includegraphics[width=170mm]{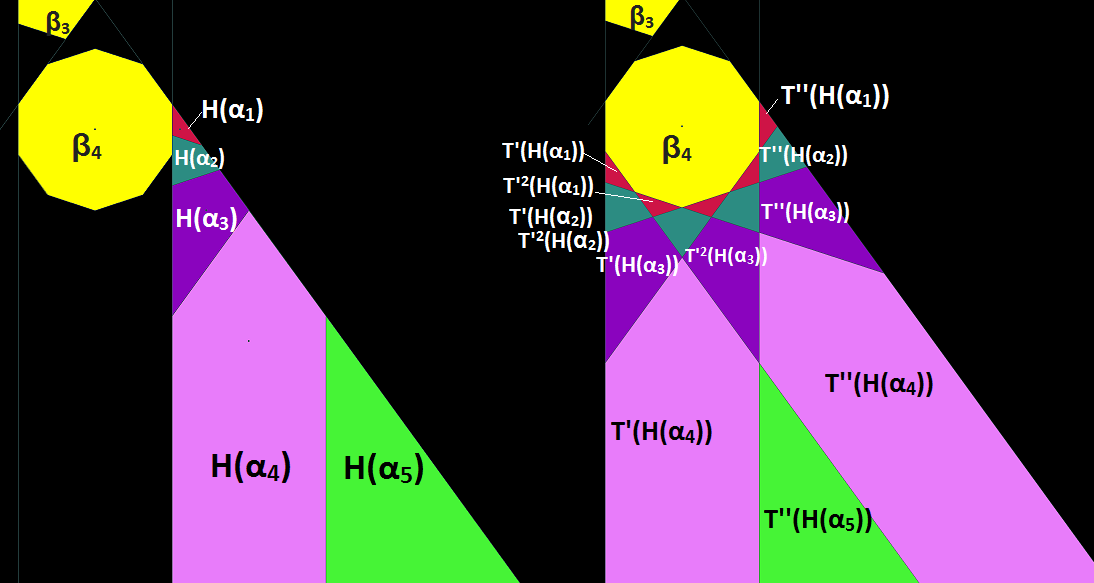}
\caption{Преобразование $T''$: траектории первого возвращения}.
\label{pic:bigSelfSimilarity-12}
\end{center}
\end{figure}

\begin{Proof}
Заметим, что $T'$ определено для всего $\alpha_5$; обозначим угол $T'(\alpha_5)$ за $\alpha'_5$. Также заметим, что правильный десятиугольник $\beta_4$ равен $\gamma$, а его центр совпадает с $O_4$ (ибо $\beta_4$ центрально-симметричен $\gamma$ относительно $Q_4$), а по лемме \ref{period1Components} точка $Q_5$, вершина угла $\alpha_5$, совпадает с вершиной $B^4_6$ многоугольника $\beta_4$. Следовательно, по лемме \ref{inducedTAction} $\alpha'_5$ есть угол, образованный лучами $A_0A_1$ и $B^4_1B^4_2$. Обозначим за  $H'$ параллельный перенос вдоль вектора $A_0B^4_0$; тогда для $\intt(V')$ верно: $H' = H \circ T'$, а $\alpha'_5 = H'(\intt(V'))$.

Пусть $V^H_{j+1} = H(V_{j+1})$, $V^{H'}_{j+1} = H'(V_{j+1})$ и $\alpha^H_{j} = H(\alpha_{j})$, $j = 1,2,3,4,5$. Тогда каждое такое $\alpha^H_j$ есть $\intt(V^H_{j+1}) \cap \alpha_5$, a $T'(\alpha^H_j)$ --- $\intt(V^{H'}_{j+1}) \cap \alpha'_5$.

Так как $\beta_4 = H'(\intt(\gamma))$, то углы $H'(V_2), \ldots, H'(V_6)$ суть углы, образованные лучами $B^4_2B^4_1$, $B^4_3B^4_2$, $\ldots$, $B^4_7B^4_6$. Более того, по лемме \ref{inducedTAction} $\forall p \in ((H'(V_2) \cup H'(V_3) \cup H'(V_4) \cup H'(V_5)) \cap \intt(V'))$: $T'(p) = R'_5(p)$, где $R'_5$ есть поворот на угол $\frac{\pi}{5}$ против часовой стрелки вокруг точки $O_4$; в частности, $R'_5(H'(V_j)) = H'(V_{j+1})$, $j = 2,3,4,5$. Так как $\alpha_5 = \intt(V^H_6)$, то тогда из всего вышесказанного следует:

\begin{itemize}
    \item $\forall j \in \{1,2,3,4,5\}$: $\forall p \in \alpha^H_j$: $T''(p)$ корректно определено, причем $T''(p) = T'^{6-j}$;
    \item $\forall j \in \{1,2,3,4,5\}$: $T''(\alpha^H_j)$ есть поворот против часовой стрелки на угол $\frac{(5-j)\pi}{5}$; в частности, $T''(\alpha^H_5)$ есть параллельный перенос на вектор $A_6A_1$;
    \item если же $p \in \alpha_5$ лежит на границе одной из $\alpha^H_j$, $j=1,2,3,4,5$, то $T'(p)$ попадет на один из лучей $B^4_3B^4_2,B^4_4B^4_3,\ldots,B^4_7B^4_6$, которые при дальнейшем применении $T'$ попадают на луч $B^4_7B^4_6$, на котором $T'$ не определено; следовательно, $T''(p)$ не определено.
\end{itemize}

Таким образом, $T''$ разделяет $\alpha_5$ на те же фигуры, что $T'$ разделяет $\intt(V')$, и применяет к ним те же повороты. Остается лишь показать, что точки $H(O_j)$, $j=1,2,3,4$, являются неподвижными. Зафиксируем $j \in \{1,2,3,4\}$. Из леммы \ref{inducedTAction} следует, что $O_j$ есть центр правильного многоугольника $\beta_j$,  вписанного в $\alpha_j$; следовательно, существует такое $r_j \in \bbr$, что открытый круг $S_j$ с центром в $O_j$ и радиусом $r_j$ вписан в $\alpha_j$. Тогда $H(S_j)$ вписан в $\alpha^H_j$; так как $\alpha_j$ есть пересечение $\intt(V')$ и $\intt(V_{j+1})$, то $H(S_j)$ вписан как в угол $\alpha_5$, так и в угол $H(V_{j+1})$, а $T'(H(S_j))$ - в угол $H''(V_{j+1})$. Так как $T'^{5-j}(\intt(H''(V_{j+1}))) = \intt(H''(V_6)) = \alpha_5$, то $T''(H(S_j)) = T'^{6-j}(H(S_j))$ есть круг, равный $H(S_j))$ и вписанная в $\alpha_5$. Следовательно, $H(O_j)$ и $T''(H(O_j))$ совпадают как центры равных вписанных в один и тот же угол кругов, QED.

\end{Proof}

Лемма \ref{selfSimilarityBig} устанавливает самоподобие для преобразования $T'$ внутри $V'$. В процессе доказательства был практически доказан важный технический факт, который мы вынесем в отдельную лемму.

\begin{Lm} \label{selfSimilarityBigStructure}
Пусть $j \in \{1,2,3,4,5\}$, и $p \in \alpha_j$. Тогда $H(T'(p)) = T'^{6-j}(H(p))$. При этом, $H(p) \in \alpha_5$, и $\forall j' \in \bbz, 0 < j' < 6-j\}$: $T'^{j'} \in \alpha_4$.
\end{Lm}

Важнейшее свойство самоподобий, похожих на рассматриваемое нами, выражает следующая лемма.

\begin{Lm} \label{typeInheritance1}
Пусть $p \in \intt(V')$, и $q = H(p)$. Тогда $p$ граничная (периодическая, апериодическая), если и только если $q$ граничная (периодическая, апериодическая).
\end{Lm}

\begin{Proof}
    Рассмотрим двусторонние последовательности $tr'(p) = \{T'^k(p) | k \in \bbz, T'^k \text{ определено} \}$, $tr'(q) = \{T'^k(q) | k \in \bbz, T'^k \text{ определено} \}$ и $tr''(q) = \{T''^k(q) | k \in \bbz, T''^k \text{ определено} \}$. Заметим, что $tr''(q)$ может быть получена из $tr'(q)$ путем удаления всех не лежащих в $\alpha_5$ точек.
    
    Из леммы \ref{selfSimilarityBig} известно, что $T'(p)$ определено, если и только если $T''(q)$ определено, причем $T''(q) = H(T(p))$; следовательно, $tr''(q) = H(tr'(p))$, и $tr''(q)$ бесконечна в две стороны, если и только если $tr'(p)$ бесконечна в две стороны. С другой стороны, если $tr'(q)$ ограничена с, например, левой стороны, т.е. для некоторого $k \in \bbz_{\leq 0}$ $T'^k(p)$ определено, а $T'^{k-1}(p)$ нет, то тогда в $tr'(q)$ встретится точка $T''^k(q)$, причем слева от нее в $tr'(q)$ будет не более четырех точек (ибо если это не так, то при последовательном применении $T'^{-1}$ к $T''^k(q)$, к точке не более четырех раз будет применен поворот по часовой стрелке на угол $\frac{\pi}{5}$ вокруг $O_4$, после чего точка окажется внутри $T'(\alpha_5)$, и следующая итерация вернет точку в $\alpha_5$); аналогичное утверждение верно и для правой стороны. Следовательно, $tr'(q)$ бесконечна в две стороны $\Leftrightarrow$ $tr''(q)$ бесконечна в две стороны $\Leftrightarrow$ $tr'(p)$ бесконечна в две стороны, откуда и следует утверждение леммы.

\end{Proof}

\end{subsection}
\begin{subsection}{Сведение проблем периодичности к ограниченному случаю}

В данном разделе мы покажем, каким образом решить проблемы периодичности и описать периодические траектории для $T'$ в углу $V'$, если мы сможем сделать это для многоугольника $Z' = A_1B^4_0B^4_9B^4_8B^4_7$. Отметим, что $Z' = \overline{\alpha_1} \cup \overline{\alpha_2} \cup \overline{\alpha_3} \cup \overline{\alpha'_4}$, где $\alpha'_4$ есть открытый четырехугольник $Q_4B^4_9B^4_8B^4_7$.

\begin{Lm}
$T'(Z') \subset Z'$.
\end{Lm}

\begin{Proof}
Заметим, что $\intt(V') \backslash \beta_4$ разбивается на две не связанные между собой связные фигуры, одной из которых является $\intt(Z')$. Тогда $T'(\alpha_j) \subset Z'$, $j = 1,2,3$ в силу того, что $O_j \in \alpha_j \subset Z'$ и $T'(O_j) = O_j$. Замечание о том, что $T'(\alpha'_4)$ есть содержащийся в $Z'$ открытый четырехугольник $B^3_1B^4_0B^4_9B^4_8$, завершает доказательство леммы.
\end{Proof}

Более того, так как $T'$ инъективно и сохраняет площадь, то верна

\begin{Lm}
$T'(Z'') \subset Z''$, где $Z'' = \overline{V' \backslash (\beta_4 \cup Z')}$.
\end{Lm}

Итак, докажем, что структура периодических/апериодических точек в $\intt(V')$ в некотором смысле порождается структурой периодических/апериодических точек в $Z' \cup \beta_4$.

\begin{Lm} \label{genPeriodicWeak}
Пусть $q \in V'$ - периодическая (апериодическая) точка. Тогда существует периодическая (апериодическая) точка $p \in Z'$, число $m \in \bbz_{\geq 0}$ и последовательность преобразований $f_1,f_2,\ldots,f_m$, т.ч. $q = f_m(f_{m-1}(\ldots f_2(f_1(p))\ldots ))$, а каждое из преобразований $f_j$, $j = 1,2,\ldots,m$ является либо $T'$, либо $H$.
\end{Lm}

\begin{Proof}
Из определения $T'$, $H$ и леммы $\ref{typeInheritance1}$ следует, что множество периодических (апериодических) точек в $V'$ инвариантно относительно $T'$ и $H$; следовательно, достаточно получить из $q$ с помощью преобразований $T'^{-1}$, $H^{-1}$ точку, лежащую в $Z'$.

Пусть $S$ есть бесконечная полоса, ограниченная отрезком $A_1B^4_1$ и лучами $A_1A_2$, $B^4_1B^4_2$. Заметим, что с точностью до граничных точек луча $A_0A_1$, $V' \backslash S = \alpha'_5 = T'(\alpha_5)$, причем для $\alpha'_5$: $T'^{-1}(r)$ есть параллельный перенос на вектор $\vec{A^1A^6}$; будем применять $T'^{-1}$ к точке $q$ до тех пор, пока $q$ не попадет в $S$. После этого, применим к $q$ $H'^{-1}$ максимальное количество раз так, чтобы точка не вышла за пределы $V'$.

По построению, теперь $q$ лежит либо в $Z' \cup \beta_4$, либо в фигуре, центрально симметричной $Z'$ относительно $O_4$. В первом случае, лемма доказана; во втором же для завершения доказательства достаточно лишь применить к $q$ преобразование $T'^{-1}$ не более трех раз так, чтобы $q$ попала в $\alpha'_5$, после чего применить еще раз $T^{-1}$; несложно показать, что в этом случае точка окажется внутри $H(Z')$, и применение $H^{-1}$ поместит точку внутрь $Z'$.

\end{Proof}

Усилим доказанную лемму.

\begin{Lm} \label{genPeriodicStrong}
Пусть $q \in V'$ - периодическая (апериодическая) точка. Тогда существует периодическая (апериодическая) точка $p \in Z'$ и числа $k, l \in \bbz_{\geq 0}$, т.ч. $q = T'^l(H^k(p))$.
\end{Lm}

\begin{Proof}
Из леммы \ref{selfSimilarityBigStructure} следует, что для любой неграничной точки $q$ верно: $H(T'(q)) = T'^j(H(q))$ для некоторого $j = j(q), j \in \bbz_+$. Для завершения доказательства достаточно применить такое равенство некоторое количество раз к равенству, полученному в лемме \ref{genPeriodicStrong}.
\end{Proof}

С помощью леммы \ref{genPeriodicStrong} можно свести ограничить область рассмотрения проблем периодичности с $V'$ до $Z'$, что мы и сделаем.

\begin{Lm} \label{reductionToZ'}
\begin{enumerate}
    \item Апериодическая точка $q \in V'$ существует, если и только если существует апериодическая точка $p \in Z'$.
    \item Периодические точки образуют внутри $V'$ множество полной меры, если и только если периодические точки образуют множество полной меры внутри $Z'$.
\end{enumerate}
\end{Lm}

\begin{Proof}
В свете леммы \ref{genPeriodicStrong}, неочевидным остается лишь утверждение о том, что если периодические точки образуют множество полной меры в $Z'$, то такие точки образуют множество полной меры и в $V'$.

Помимо периодических точек, в $V'$ существуют лишь точки граничные и апериодические. Так как множество граничных точек внутри $Z'$, $V'$ и в целом в $\outtable$ имеет меру нуль по лемме \ref{zeroMeasureOfBoundaryPoint}), то остается лишь доказать, что если апериодические точки внутри $Z'$ имеют меру нуль, то такой мерой обладают и апериодические точки внутри всего угла $V'$.

По лемме \ref{genPeriodicStrong}, каждую апериодическую точку $q \in V'$ можно представить в виде $q = T'^l(H^k(p))$, $p = p(q)$ - апериодическая точка, $k = k(q), l = l(q), k,l \in \bbz_{\geq 0}$. Заметим, что можно рассмотреть процесс вычисления $(H^k \circ T'^l)(p)$ как последовательное применение к точке $p$ $k+l$ движений плоскости (всей), каждое из которых является либо параллельным переносом на один из двух векторов ($\vec{P_1Q_5}$ либо $\vec{A_6A_1}$), либо одним из поворотов против часовой стрелки на угол $\frac{(5-j)\pi}{5}$ вокруг точки $O_j$, $j=1,2,3,4$ (см. определение $H$ и лемму \ref{inducedTAction}). Т.к. композиция движений есть движение, то мы можем сопоставить каждой апериодической точке $q$ движение из счетного множества $\Phi$ движений вышеописанного типа.

Пусть $V'_{\infty}$, $Z'_{\infty}$ суть множества апериодических точек внутри $V'$ и $Z'$ соответственно. Тогда по построению и \ref{genPeriodicStrong}: $V'_{\infty} \subset \bigcup\limits_{f \in \Phi} f(Z'_{\infty})$. Так как $f$ есть изометрия, а $Z'_{\infty}$ имеет меру нуль, то и $f(Z'_{\infty})$ имеет меру нуль; следовательно, $V'_{\infty}$ есть подмножество счетного объединения множеств нулевой меры и, как следствие, само обладает мерой нуль, QED.
\end{Proof}

Отметим, что в терминах доказательства предыдущей леммы, если бы мера $Z'_{\infty}$ была положительной, то $V'_{\infty}$ была бы бесконечной, ибо $V'_{\infty}$ содержит в себе как подмножество бесконечное число непересекающихся множеств  $\bigcup\limits_{k \in \bbz_{\geq 0}}H^k(Z'_{\infty})$.

\end{subsection}

\begin{subsection}{Сведение нахождения периодов к ограниченному случаю: введение подстановки}

Свести нахождение множества периодов для $V'$ к нахождению множества периодов для $Z' \cup \beta_4$ нам поможет символическая динамика.

Введем несколько определений, базируясь на \cite{PF02}. Пусть $A, B$ --- конечные множества, называемые {\it алфавитами}.

\begin{Def}
Конечным словом длины $l \in \bbz_{\geq 0}$ над алфавитом $A$ назовем последовательность $u_0u_1\ldots u_{l-1}$, т.ч. $u_i \in A$, $i = 0, 1, \ldots, l-1$. Множество всех конечных слов над $A$ обозначим как $A^*$.
\end{Def}

\begin{Def}
Бесконечным в две стороны словом над алфавитом $A$ назовем последовательность $(u_i)_{i \in \bbz} = \ldots u_{-2}u_{-1}u_0u_1u_2 \ldots $, т.ч. $\forall i \in \bbz: u_i \in A$. Множество всех бесконечных в две стороны слов над $A$ обозначим как $A^{\bbz}$.
\end{Def}

Аналогичным образом можно ввести и бесконечные в одну сторону слова; однако для наших целей ограничимся бесконечными в две стороны словами. Отметим, что такими словами являются коды периодических и апериодических точек; а именно:

\begin{itemize}
    \item если $p \in \outtable$ --- неграничная точка, то $\rho(p) \in \{0,1,\ldots,9\}^{\bbz}$;
    \item если при этом $p \in V'$, то $\rho'(p) \in \{0,1,\ldots,4\}^{\bbz}$.
\end{itemize}

Аналогичным разделу \ref{basicsOfPolygonalOuterBilliard} образом для конечного или бесконечного слова $U$ будем употреблять обозначения $U[i]$ и/или $U[l, r]$, если целые индексы $i, l, r$ не выводят нас за пределы слова $U$ и $l \leq r$.

\begin{Def}
Пусть $U,V \in A^*$, $U = u_0u_1\ldots u_{l-1}$, $V = v_0v_1\ldots, v_{m-1}$. Тогда конкатенацией слов $U, V$ назовем слово $UV = u_0u_1\ldots u_{l-1}v_0v_1\ldots v_{m-1}$.
\end{Def}

\begin{Def}
Пусть $\forall i \in \bbz: U_i \in A^* \backslash \{\epsilon\}$ --- непустое слово длины $l_i$. Тогда конкатенацией слов $\ldots, U_{-2}, U_{-1}, U_0, U_1, U_2, \ldots$ назовем бесконечное в две стороны слово $U = \ldots U_{-2}U_{-1}U_0U_1U_2 \ldots$, устроенное следующим образом. Пусть $m_i$, $i \in \bbz$ --- такая бесконечная в две стороны последовательность целых чисел, что:

\begin{itemize}
    \item $m_0 = 0$;
    \item $\forall i \in \bbz_+$: $m_i = m_{i-1} + l_{i-1}$;
    \item $\forall i \in \bbz_-$: $m_i = m_{i+1} - l_i$.
\end{itemize}

Тогда требуемая конкатенация есть такое слово $U \in A^{\bbz}$, что $\forall i \in \bbz$: $U[m_i, m_{i+1}-1] = U_i$.
\end{Def}

Введем теперь понятие подстановки, играющее ключевую роль в нашем исследовании.

\begin{Def}
Пусть $\sigma: A \rightarrow B^* \backslash \{\epsilon\}$ --- произвольная функция, где $\epsilon$ --- пустое слово. Расширим ее до $\sigma: A^* \cup A^{\bbz} \rightarrow B^* \cup B^{\bbz}$ c помощью следующих правил:

\begin{itemize}
    \item $\sigma(\epsilon) = \epsilon$;
    \item если $W = u_0u_1\ldots u_{l-1} \in A^*$, то $\sigma(W) = \sigma(u_0)\sigma(u_1)\ldots \sigma(u_{l-1}) \in B^*$;
    \item если $W = \ldots u_{-2}u_{-1}u_0u_1u_2\ldots \in A^{\bbz}$, то $\sigma(W) = \ldots \sigma(u_{-2})\sigma(u_{-1})\sigma(u_0)\sigma(u_1)\sigma(u_2) \ldots   $.
\end{itemize}

Таким образом устроенную функцию $\sigma$ будем называть подстановкой.
\end{Def}

Для удобства будем иногда говорить, что $\sigma$ определена на символах алфавита $A$, т.е. $A \rightarrow B^* \backslash \{\epsilon\}$, подразумевая ее определение на $A^* \cup A^{\bbz}$.

Вернемся к индуцированному внешнему биллиарду.

\begin{Def}
Пусть $p \in \outtable (V')$ --- периодическая точка с периодом $l (l')$ относительно $T (T')$. Тогда кодом периода $\rho_{per}(p) (\rho'_{per}(p))$ назовем слово $\rho(p)[0, l-1] (\rho'(p)[0, l'-1])$. \end{Def}

Пусть $\psi$ есть подстановка, определенная над алфавитом $\{1,2,3,4,5\}$ таким образом, что:

\begin{itemize}
    \item $\psi(1) = 54444$;
    \item $\psi(2) = 5444$;
    \item $\psi(3) = 544$;
    \item $\psi(4) = 54$;
    \item $\psi(5) = 5$.
\end{itemize}

Тогда из структуры преобразования $T''$ первого возвращения на $\alpha_5$ (см. лемму \ref{selfSimilarityBig}) напрямую следует

\begin{Lm} \label{selfSimilarityPeriodCodeSubstitution}
Пусть $p \in V'$ - периодическая или апериодическая точка. Тогда:

\begin{itemize}
    \item $\rho'(H(p)) = \psi(\rho'(p))$;
    \item $\rho'_{per}(H(p)) = \psi(\rho'_{per}(p))$.
\end{itemize}
\end{Lm}

Лемма \ref{selfSimilarityPeriodCodeSubstitution}, вкупе с леммой \ref{genPeriodicStrong}, позволяет найти все возможные коды периодов периодических точек множества $V'$, если таковые коды найдены для множества $T'$ с помощью следующей, очевидно следующей из указанных двух, леммы.

\begin{Lm} \label{fromZ'PeriodCodesToV'PeriodCodes}
Пусть $P_{Z'\cup \beta_4} (P_{V'}) = \{\rho'_{per}(p)\ |\ p \in Z'\cup \beta_4 (V'),\ p - \text{ периодическая компонента}\}$. Тогда $P_{V'} = \{S^l(\sigma^k(W))\ |\ W \in P_{Z' \cup \beta_4}, k,l \in \bbz_{\geq 0} \}$, где $S$ есть циклический сдвиг слова (т.е. $S(u_0u_1\ldots u_{m-1}) = u_1u_2\ldots u_{m-1}u_0$ для произвольного слова $u_0u_1\ldots u_{m-1}$ над произвольным конечным алфавитом).
\end{Lm}

Для удобства поиска самих периодов, введем еще два понятия, используемые в \cite{PF02}.

\begin{Def} \label{def:abelization}
Пусть $A = \{a_1, a_2, \ldots, a_d\}$ - конечный алфавит размера $d$. Тогда каноническим гомоморфизмом, или гомоморфизмом абелизации, назовем такое преобразование $c: A^* \rightarrow \bbz^d$, что для произвольного слова $W \in A^*$: $c(W)$ есть столбец, $i$-я координата которого есть количество раз, которое символ $a_i$ встречается в $W$, $i = 1,2,\ldots,d$.
\end{Def}

Отметим, что если $p \in V' (\outtable)$ есть периодическая точка с кодом периода $w' = \rho'_{per}(p) (w = \rho_{per}(p))$, то период точки $p$ относительно $T'(T)$ есть $(1,1,1,1,1)*c(w')$ ($(1,1,1,1,1,1,1,1,1,1)*c(w)$).

\begin{Def}
Пусть $A = \{a_1, a_2, \ldots, a_d\}$, $B = \{b_1, b_2, \ldots, b_m\}$ --- конечные алфавиты, а $\sigma: A \rightarrow B \backslash \{\epsilon\}$ --- произвольная подстановка. Тогда матрицей подстановки $M_{\sigma}$ назовем матрицу $c(\sigma(a_1)), c(\sigma(a_2)), \\ \ldots, c(\sigma(a_d)))$. Другими словами, $M_{\sigma} = ||c_{ij}||_{m \times d}$, где $c_{ij}$ есть количество раз, которое символ $b_i$ встречается в слове $\sigma(a_j)$.
\end{Def}

Например, для вышеопределенной подстановки $\psi$:

{ \small
$
     M_{\psi} = \begin{pmatrix}
0 & 0 & 0 & 0 & 0 \\
0 & 0 & 0 & 0 & 0 \\
0 & 0 & 0 & 0 & 0 \\
4 & 3 & 2 & 1 & 0 \\
1 & 1 & 1 & 1 & 1
\end{pmatrix}.
$

} % \Small

В терминах гомоморфизма абелизации и матрицы инцидентности, леммы \ref{fromIPStructureToPeriod} и \ref{fromZ'PeriodCodesToV'PeriodCodes} могут быть переформулированы следующим образом.

\begin{Lm} \label{fromIPStructureToPeriodAbel}
Пусть $p \in V'$ - периодическая точка или компонента. Пусть $c(p) = c(\rho_{per}(p))$, а $c'(p) = c(\rho'_{per}(p))$. Тогда $per(p) = \frac{10*((1,1,1,1,1)*c'(p))}{\rgcd(10, (1,2,3,4,5))*c'(p)}$.
\end{Lm}

\begin{Lm} \label{fromZ'PeriodCodesToV'PeriodCodesAbel}
Пусть $C_{Z'\cup \beta_4} (C_{V'}) = \{c(\rho'_{per}(p))\ |\ p \in Z'(V'), p - \text{ периодическая компонента}\}$. Тогда $C_{V'} = \{M_{\psi}^k(w)~|~w \in C_{Z'\cup \beta_4}, k \in \bbz_{\geq 0} \}$. 
\end{Lm}

\end{subsection}

\begin{subsection}{$Z' \cup \beta_4$: описание структуры}

Леммы \ref{reductionToZ'}, \ref{fromIPStructureToPeriodAbel} вкупе с фактом, что $\beta_4$ есть периодическая компонента с кодом периода $\rho'_{per}(\beta_4) = 4$, позволяет нам сосредоточить все внимание на $Z'$.

Начнем со следующего замечания. Пусть $Z'_1$ есть треугольник $P_1P_3Q_3$, а $Z'_2$ есть семиугольник $B^3_3B^3_2B^3_1B^4_0B^4_9B^4_8B^4_7$, лежащий внутри $Z'$ <<между>>  $\beta_3$ и $\beta_4$.

\begin{Lm}
$Z'$ можно разбить на многоугольники-фигуры $Z'_1 = P_1P_3Q_3$, $\beta_3$ и $Z'_2 = B^3_3B^3_2B^3_1B^4_0B^4_9B^4_8B^4_7$, каждая из которых $T'$-инвариантна.
\end{Lm}

\begin{Proof}
Утверждение про разбиение очевидно из рис. \ref{pic:InducedT-InvariantComponents-12}. Чтобы доказать инвариантность, достаточно заметить, что:
\begin{itemize}
    \item $Z'_1$ и $Z'_2$ не имеют общих точек (даже на границе);
    \item $T'$ делит $Z'_1$ на две связные фигуры $\overline{\alpha_1}$ и $\overline{\alpha_2}$, для каждой из которых $T'$ есть поворот вокруг содержащейся в фигуре неподвижной точки (см. лемму \ref{inducedTAction});
    \item $T'$ делит $Z'_2$ на две связные фигуры, являющиеся $\overline{\alpha'_4}$ и треугольником $B^3_1B^4_0B^3_2$, причем $T'(\alpha'_4)$ есть открытый четырехугольник $B^3_1B^4_0B^4_9B^4_8 \subset Z'_2$, а $T'(\intt(B^3_1B^4_0B^3_2))$ есть открытый треугольник $B^3_2B^4_7B^3_3 \subset Z'_2$ (по той же лемме \ref{inducedTAction})).
\end{itemize}
\end{Proof}

$\beta_3$ есть периодическая компонента с кодом периода $\rho'_{per}(\beta_3) = 3$; оставшиеся периодические и апериодические точки, компоненты и орбиты можно искать независимо в $Z'_1$ и $Z'_2$.

Рассмотрим преобразование $T'$ на фигуре $Z'_2$. Пусть $X'_2$ есть открытый четырехугольник $B^3_3B^4_9B^4_8B^4_7$. Пусть $T'_2: X'_2 \rightarrow X'_2$ есть преобразование первого возвращения относительно $T'$.

\begin{figure}[h!]
\begin{center}
\includegraphics[width=170mm]{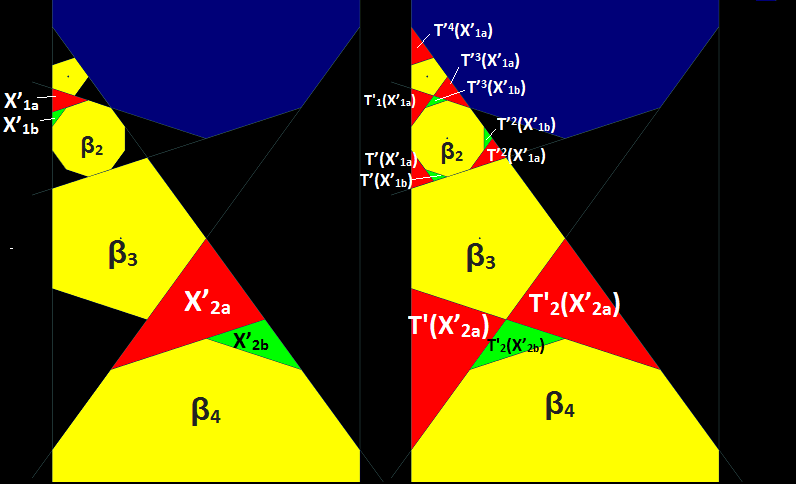}
\caption{Преобразования $T'_1, T'_2$: траектории первого возвращения}.
\label{pic:X12-12}
\end{center}
\end{figure}

Из рис. \ref{pic:X12-12} очевидно следует доказательство следующей леммы.

\begin{Lm} \label{X'_2Structure}
\begin{enumerate}
    \item $T'_2$ делит $X'_2$ на два открытых треугольника $X'_{2a}$ и $X'_{2b}$ (точный смысл обозначений будет ясен позже) лучом $B^4_9B^4_8$, причем $X'_{2a}$ имеет большую площадь; для точек луча $B^4_9B^4_8$ внутри $X'_2$ $T'_2$ не определено; 
    \item $T'_2(X'_{2a})$ есть поворот на угол $\frac{3\pi}{5}$ против часовой стрелки вокруг точки $O'_{2u}$ пересечения биссектрис углов четырехугольника $X'_2$, причем $T'_2(X'_{2a}) = T'^2(X'_{2a})$;
    \item $T'_2(X'_{2b})$ есть поворот на угол $\frac{\pi}{5}$ против часовой стрелки вокруг точки $O_4$, причем $T'_2(X'_{2b}) = T'(X'_{2b})$;
    \item $\overline{Z'_2} = \overline{X'_{2a}} \cup \overline{T'(X'_{2a})} \cup \overline{X'_{2b}}$, причем границы всех упомянутых многоугольников состоят из граничных точек.
\end{enumerate}
\end{Lm}

Рассмотрим теперь $Z'_1$, а внутри него - открытый четырехугольник $X'_1$, равный $B^1_1B^2_0B^2_9B^2_8$. Пусть $T'_1: X'_1 \rightarrow X'_1$ - преобразование первого возвращения относительно $T'$. Из все того же рис. \ref{pic:X12-12} очевидно следует доказательство следующей, аналогичной предыдущей, леммы.

\begin{Lm} \label{X'_1Structure}
\begin{enumerate}
    \item $T'_1$ делит $X'_1$ на два открытых треугольника $X'_{1a}$ и $X'_{1b}$ лучом $B^2_8B^2_9$, причем $X'_{1a}$ имеет большую площадь; для точек луча $B^2_8B^2_9$ внутри $X'_1$ $T'_1$ не определено;
    \item $T'_1(X'_{1a})$ есть поворот на угол $\frac{3\pi}{5}$ по часовой стрелке вокруг точки $O'_{1u}$ пересечения биссектрис углов четырехугольника $X'_1$, причем $T'_1(X'_{1a}) = T'^5(X'_{1a})$;
    \item $T'_1(X'_{1b})$ есть поворот на угол $\frac{\pi}{5}$ по часовой стрелке вокруг точки $O_2$, причем $T'_1(X'_{1b}) = T'^3(X'_{2b})$;
    \item $\overline{Z'_1} = \bigcup\limits_{j=0}^4 \overline{T'^j(X'_{1a})} \cup \bigcup\limits_{j=0}^2\overline{T'(X'_{2b})} \cup \beta_1 \cup \beta_2$, причем границы всех упомянутых многоугольников состоят из граничных точек.
\end{enumerate}
\end{Lm}

Таким образом, любая периодическая или апериодическая траектория внутри $Z'_1 \backslash (\beta_1 \cup \beta_2)$ или $Z'_2$ обязательно проходит через $X'_1$ или $X'_2$, причем на этой траектории не может существовать больше четырех идущих подряд точек, не лежащих внутри $X'_1$ или $X'_2$. Отсюда и из того факта, что $\beta_1, \beta_2$ суть периодические компоненты, очевидно следует

\begin{Lm} \label{reductionToX12}
\begin{enumerate}
    \item Пусть $p \in X'_j$, $j=1,2$. Тогда $p$ граничная (периодическая, апериодическая) относительно $T'$, если и только если $p$ граничная (периодическая, апериодическая) относительно $T'_j$;
    \item Апериодическая точка внутри $Z'$ существует, если и только если такая точка существует внутри $X'_1$ или $X'_2$;
    \item Периодические точки образуют в $Z'$ множество полной меры, если и только если эти точки образуют множество полной меры в $X'_1$ и $X'_2$.
\end{enumerate}
\end{Lm}

\end{subsection}
\begin{subsection}{Преобразования первого возвращения и динамическая система $(X, f)$}

Согласно лемме \ref{reductionToX12}, изучение проблем периодичности сведено к преобразованиям $T'_1$ и $T'_2$ на фигурах $X'_1$ и $X'_2$. Заметим, что преобразования $T'_1$ и $T'_2$ <<одинаковы>> с точностью до симметрии, что дает возможность изучить их одновременно. Для этого рассмотрим произвольный открытый четырехугольник $X = ABCD \subset \bbrr$, вершины которого перенумерованы против часовой стрелки, подобный $X'_2 = B^3_3B^4_9B^4_8B^4_7$ и $X'_1$. Введем также такие аффинные преобразования плоскости $\Delta_1, \Delta_2$, что $\Delta_1(X) = X'_1$, $\Delta_2(X) = X'_2$, причем:

\begin{itemize}
    \item $\Delta_1(A) = B^1_1$, $\Delta_1(B) = B^2_8$, $\Delta_1(C) = B^2_9$, $\Delta_1(D) = B^2_0$;
    \item $\Delta_2(A) = B^3_3$, $\Delta_2(B) = B^4_9$, $\Delta_2(C) = B^4_8$, $\Delta_2(D) = B^4_7$.
\end{itemize}

Отметим, что $\Delta_1$ меняет ориентацию, а $\Delta_2$ - нет.

Пусть теперь $E$ есть пересечение прямых $BC$ и $AD$, а $F$ --- пересечение прямых $AB$ и $CD$. Введем кусочно-аффинное (а точнее, <<кусочно-вращающее>>) преобразование $f: X \rightarrow X$, т.ч.:

\begin{itemize}
    \item если $p \in \intt(ABE)$, то $f(p)$ есть поворот плоскости, переводящий треугольник $ABE$ в треугольник $FDA$;
    \item если $p \in \intt(CED)$, то $f(p)$ есть поворот плоскости, переводящий треугольник $CED$ в треугольник $BFC$;
    \item если $p \in \partial(ABE) \cup \partial(CED)$, то $f(p)$ не определено.
\end{itemize}

Фигура $X$ и преобразование $f$ изображены на рис. \ref{pic:Xf-12}.

\begin{figure}[h!]
\begin{center}
\includegraphics[width=170mm]{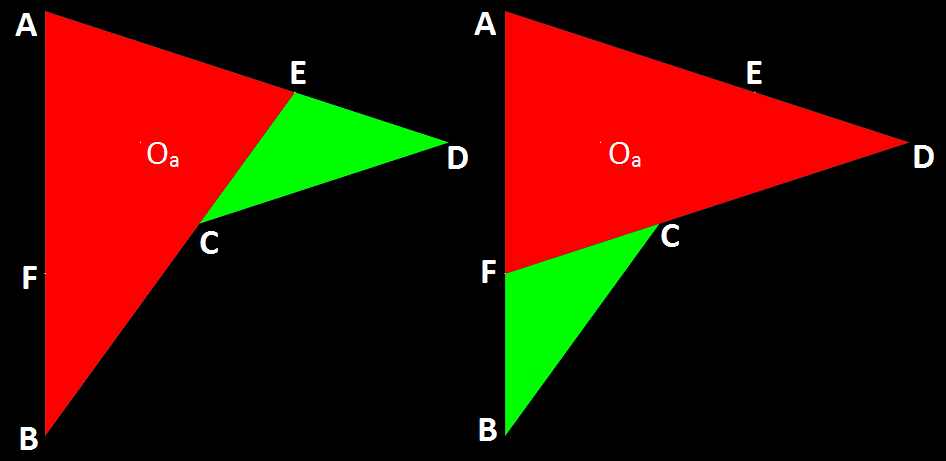}
\caption{Фигура $X$ и преобразованиe $f$: траектории первого возвращения}.
\label{pic:Xf-12}
\end{center}
\end{figure}

Очевидным образом вводятся также понятия граничности, периодичности и апериодичности и периодической компоненты внутри $X$ относительно $f$. Более того, для неграничной точки $p \in X$ введем также {\it код} $\rho_X(p) = \ldots w_{-2}w_{-1}w_0w_1w_2 \ldots$, т.ч. $\forall j \in \bbz$, $w_j \in \{a, b\}$, причем
\begin{itemize}
    \item $a$, если $f^j \in \intt(ABE)$;
    \item $b$, если $f^j \in \intt(CED)$.
\end{itemize}

Здесь $\{a, b\}$ есть двухбуквенный алфавит.

Из структуры $X,f,X'_1,T'_1,X'_2,T'_2$ и лемм \ref{X'_1Structure}, \ref{X'_2Structure} следует лемма, иллюстрирующая окончательное сведение задачи к фигуре Х.

\begin{Lm} \label{reductionToX}
Пусть $p \in X$, $q_1 = \Delta_1(p)$, $q_2 = \Delta_2(p)$. Тогда:
\begin{itemize}
    \item $p$ граничная (периодическая, апериодическая) относительно $f$, если и только если $q_j$, $j=1,2$, граничные (периодические, апериодические) относительно $T'$;
    \item Для каждой из проблем периодичности, ответ для $X'_1$ и $X'_2$ относительно $T'$ положительный, если и только если ответ для $X$ относительно $f$ положительный;
    \item если $p$ периодическая или апериодическая, то:
        \begin{itemize}
            \item $\rho'(q_1) = \phi_1 (\rho_X(p))$, где $\phi_1: \{a,b\} \rightarrow \{1,2,3,4,5\}^*$ - подстановка, т.ч. $\phi_1(a) = 22211$, $\phi_1(b) = 222$;
            \item $\rho'(q_2) = \phi_2 (\rho_X(p))$, где $\phi_2: \{a,b\} \rightarrow \{1,2,3,4,5\}^*$ - подстановка, т.ч. $\phi_2(a) = 43$, $\phi_2(b) = 4$.
        \end{itemize}
    \item если $p$ периодическая, то $\rho'_{per}(q_j) = \phi_j(\rho_{Xper}(p))$, $j=1,2$;
    \item пусть $C_X = \{c(\rho_{Xper}(p)) | p \subset X, p \text{ - периодическая компонента}\}$. Тогда $C_{Z'\cup \beta_4} = \{M_{\phi_j}(w) | w \in C_X, j = 1,2\} \cup \{\begin{pmatrix} 1 & 0  & 0 & 0 & 0 \end{pmatrix}^T, \begin{pmatrix} 0 & 1 & 0 & 0 & 0 \end{pmatrix}^T, \\ \begin{pmatrix} 0 & 0 & 1 & 0 & 0 \end{pmatrix}^T, \begin{pmatrix} 0 & 0 & 0 & 1 & 0 \end{pmatrix}^T \}$. 
\end{itemize}
\end{Lm}

\end{subsection}
\begin{subsection}{Фигура X: базовые периодические компоненты}

В следующих нескольких разделах будут исследоваться введенные в предыдущем разделе фигура $X$ и преобразование $f$ на ней. Динамическая система $(X, f)$ возникла в работе С.Табачникова \cite{Tab93}, ибо в случае внешнего биллиарда вне правильного пятиугольника динамическая система, аналогичная $(Z', T')$ для случая десятиугольника, изоморфна $(X, f)$; как следствие, $(X, f)$ была тщательно изучена. Тем не менее, для полноты картины мы проведем исследование для системы $(X, f)$ <<с нуля>>.

Пусть $O_a$ есть точка пересечения биссектрис углов $ABE$ и $FDA$. В силу симметрии очевидно, что $O$ лежит также и на отрезке $AC$, также являющемся частью биссектрисс углов $BAD$ и $ECF$. Отсюда следует, что существует открытый десятиугольник $\omega_a$ с центром в $O_a$, вписанный в четырехугольник $AFCE$. Более того, так как угол $FCE$ по построению равен $\frac{3\pi}{5}$, углу правильного десятиугольника, то одна из вершин $\omega_a$ совпадает с вершиной этого угла, т.е. точкой $C$. Обозначим вершины $\omega_a$ $W_0, W_1, \ldots, W_9$, нумеруя вершины против часовой стрелки таким образом, что $W_3 = C$. Тогда точки $W_2$, $W_4$ лежат соответственно на отрезках $FC$ и $CE$, а отрезки $W_6W_7$ и $W_9W_0$ - на отрезках $DA$ и $AB$.

Важность десятиугольника $\omega_a$ иллюстрирует следующая лемма.

\begin{Lm}
\begin{enumerate}
    \item $f(O_a) = O_a$;
    \item $f(\omega_a) = \omega_a$;
    \item $\partial \omega_a$ состоит из граничных точек;
    \item $\omega_a$ есть периодическая компонента относительно $f$ с кодом периода $\rho_{Xper}(\omega_a) = a$.
\end{enumerate}
\end{Lm}

\begin{Proof}
Из определения $f$ следует, что $f$ переводит биссектрису угла $ABE$ в биссектрису угла $FDA$, а точку $O_a$ - в точку $f(O_a)$, лежащую на биссектрисе угла $FDA$, причем $|AO_a| = |Cf(O_a)|$; из симметрии $X$ и определения $O_a$ следует, что $O_a = f(O_a)$, что доказывает п.1.

Так как $O_a = f(O_a)$, то для треугольника $ABE$ $f$ есть поворот на угол $\frac{3\pi}{5}$. Тогда п.2 леммы верен уже потому, что $\omega_a$ есть правильный десятиугольник с центром в $O_a$, целиком лежащий в треугольника $ABE$. Более того, при таком повороте стороны $\omega_a$ переходят друг в друга, перескакивая против часовой стрелки через две стороны на третью. Так как числа 3 и 10 взаимно просты, то при последовательном применении поворота каждая из сторон $\omega_a$ побывает каждой. Так как одна из сторон $\omega_a$ лежит на отрезке $CF$, для которой $f$ не определено, то все стороны $\omega_a$ состоят из граничных точек, и п.3 доказан. П.4 следует напрямую из пп.2,3. 
\end{Proof}

Так как $T(\intt(CED)) \subset \intt(ABE)$, то апериодических компонент с кодом периода $b$ для $(X, f)$ нет. Однако рассмотрим точку $O_{ab}$, являющуюся точкой пересечения биссектрис треугольника $CED$. Пусть $G$ есть точка пересечения отрезка $CE$ и луча $W_6W_5$, а $H$ = точка пересечения лучей $W_0W_1$ и отрезка $FC$. В силу симметрии относительно биссектрисы $DO_a$ угла $CDE$, получаем, что $O_{ab}$ лежит также и на биссектрисе угла $GB_6D$.

По определению, $f(\intt(CED)) = \intt(BFC)$; очевидно, что $O_{ba}$ перейдет в точку $O_{ab}$ пересечения биссектрис треугольников $BW_0H$ и $BFC$. С другой стороны, $f(\intt(BFC)) = \intt(DGW_6)$, причем точкой пересечения биссектрис треугольника $DGW_6$ является $O_{ba}$. Следовательно, $f^2(O_{ba}) = O_{ba}$.

Более того, так как $CED$ есть равнобедренный треугольник с углом $\frac{3\pi}{5}$ при вершине $E$, то существует открытый правильный пятиугольник $\omega_{ba}$ с центром в $O_{ba}$, вписанный в треугольник $CED$. Так как $f^2(CED)$ есть поворот на угол $\frac{4\pi}{5}$ вокруг точки $O_{ba}$, то $f^2(\omega_{ba}) = \omega_{ba}$, а стороны $\omega_{ba}$ при последовательном применении $f^2$ <<перескакивают>> через одну и рано или поздно попадут на отрезок $CD$ либо $DE$. Аналогичными свойствами обладает и пятиугольник $\omega_{ab} = f(\omega_{ba})$.

Подведем итог предыдущих абзацев следующей, вытекающей из вышесказанного, леммой.

\begin{Lm}
\begin{enumerate}
    \item $f(O_{ba}) = O_{ab}$, $f(O_{ab}) = O_{ba}$;
    \item $f(\omega_{ba}) = \omega_{ab}$, $f(\omega_{ab}) = \omega_{ba}$;
    \item $\omega_{ba} (\omega_{ab})$ есть периодическая компонента относительно $f$ с кодом периода $ba$ ($ab$).
\end{enumerate}
\end{Lm}

\end{subsection}
\begin{subsection}{Фигура Х: самоподобие 2}

Пусть $\Gamma: \bbrr \rightarrow \bbrr$ - такое аффинное преобразование, что $\Gamma(A) = A$, $\Gamma(B) = W_7$, $\Gamma(C) = W_8$, $\Gamma(D) = W_9$; другими словами, $\Gamma$ есть композиция сжатия в $\frac{|AB|}{AW_9}$ раз с центром в точке $A$ и осевой симметрии относительно биссектрисы $AC$ угла $BAD$. Пусть $X_{\Gamma} = \Gamma(X)$. Пусть $f_{\Gamma}: X_{\Gamma} \rightarrow X_{\Gamma}$ есть преобразование первого возвращения на $X_{\Gamma}$ относительно $f$ (см. определение \ref{bigFirstReturnMap} аналогичного преобразование первого возвращения на $\alpha_5$ относительно $T'$).

Самоподобие системы $(X, f)$ устанавливает следующая лемма.

\begin{Lm} \label{selfSimilarityX}
Пусть $p \in X$, и $q = \Gamma(p)$. Тогда:
\begin{itemize}
    \item $f_{\Gamma}(q)$ определено, если и только если $f(p)$ определено;
    \item в случае, если $f(p)$ определено, $f_{\Gamma}(q) = \Gamma(f(p))$;
    \item если $p \in \intt(ABE)$, то $f_{\Gamma}(q) = f^7(q)$, а если $p \in \intt(CED)$, то $f_{\Gamma(q)} f^3(q)$. 
\end{itemize}
\end{Lm}

Доказательство ясно из рисунка \ref{pic:Xf-SelfSimilarity}, на котором изображены открытые треугольники $f^i(\Gamma(\intt(ABE)))$, $i = 0, 1, \ldots, 6$, а также открытые треугольники $f^j(\Gamma(\intt(CED)))$, $j = 0, 1, 2$. 

\begin{figure}[h!]
\begin{center}
\includegraphics[width=170mm]{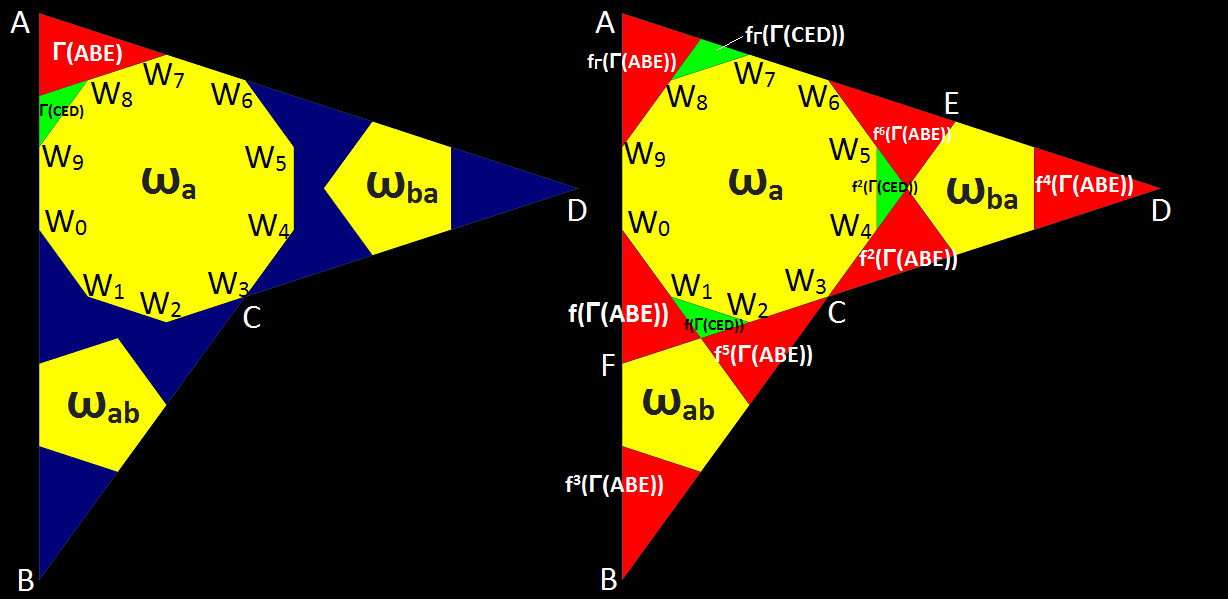}
\caption{Фигура $X$ и преобразованиe $f_\Gamma$: траектории первого возвращения}.
\label{pic:Xf-SelfSimilarity}
\end{center}
\end{figure}

Аналогично лемме \ref{typeInheritance1} можно доказать и следующую лемму.

\begin{Lm} \label{typeInheritance2}
Пусть $p \in X$, $q = \Gamma(p)$. Тогда $p$ граничная (периодическая, апериодическая), если и только если $q$ граничная (периодическая, апериодическая).
\end{Lm}

Введем теперь еще несколько определений.

\begin{Def}
Пусть $q \in X$. Тогда рангом $rk(q)$ точки $q$ назовем максимальное целое неотрицательное число $k$, т.ч $\Gamma^{-k}(q) \in X$.  
\end{Def}

\begin{Def}
Пусть $q \in X$ - периодическая точка. Тогда рангом орбиты $rko(q)$ назовем число $\max\limits_{j=0}^{\rho_{Xper} - 1} rk(f^j(q))$.
\end{Def}

\begin{Lm} \label{fromRankToRank}
Пусть $q \in X$ - периодическая точка, и $rk(q) = rko(q) = k > 0$. Тогда $rk(p) = rko(p) = k - 1$, где $p = \Gamma^{-1}(q)$.
\end{Lm}

\begin{Proof}
По лемме \ref{typeInheritance2} $p$ также периодическая; пусть ее период относительно $f$ равен $l$. Тогда по лемме \ref{selfSimilarityX}: $\forall j, 0 \leq j < l: f'^j(q) = \Gamma(f^j(p)))$, причем $f'^l(q) = q$. Согласно определениям $f'$ и ранга, среди точек орбиты точки $q$ относительно $f$ только точки $\{f'^j(q), 0 \leq j < l\}$ обладают отличным от нуля рангом, причем $\forall j \in [0, l)$: $rk(f'^j(q)) = 1 + rk(f^j(p))$. Следовательно, $rko(p) = rk(p) = k-1$, QED.    
\end{Proof}

Из леммы \ref{fromRankToRank} следует следующая, аналогичная лемме \ref{genPeriodicStrong}, лемма. 

\begin{Lm} \label{genPeriodicStrong2}
Любая периодическая точка (компонента) $q \in (\subset) X$ может быть представлена в виде $q = f^l(\Gamma^k(p))$, где $p \in (\subset) X \backslash \Gamma(X)$, $p$ - периодическая точка ранга 0 ($p \in \{\omega_a, \omega_{ba}\}$, $k, l \in \bbz_{\geq 0}$. 
\end{Lm}

\begin{Proof}
<<Точечная>> часть очевидно следует из леммы \ref{fromRankToRank}.
<<Компонентная>> же часть леммы \ref{genPeriodicStrong2} следует из того, что граница $(\partial \Gamma(X)) \cap X$ состоит лишь из граничных точек. Из этого же утверждения следует, что понятия ранга и ранга орбиты могут быть естественным образом обобщены и на периодические компоненты внутри $X$. Из рис. \ref{pic:Xf-SelfSimilarity} следует, что периодическими компонентами ранга 0 являются правильные многоугольники $\omega_a, \omega_{ba}$ и $\omega_{ab} = f(\omega_{ba})$, ибо все периодические точки $X$, не принадлежащие $\omega_a, \omega_{ba}, \omega_{ab}$ лежат на траекториях возвращения $f$ в $\Gamma(X)$.
\end{Proof}
\end{subsection}
\begin{subsection}{Доказательство теоремы 1}
В данном разделе мы установим существование апериодической точки для системы $(X, f)$; как следует из лемм \ref{reductionToX}, \ref{reductionToZ'}, \ref{reductionToT'}, этого будет достаточно для доказательства теоремы \ref{MainTreorem}.

Пусть $\Gamma_1 = \Gamma \circ f^{-1} \circ \Gamma \circ f $; заметим, что последовательность $X, \Gamma_1(X), \Gamma_1^2(X), \Gamma_1^3(X), \ldots $ есть последовательность вложенных друг в друга четырехугольников, которая сходится к некоторой точке $p_{\inf}$, лежащей строго внутри каждого из $\Gamma_1^j(X), j \in \bbz_{\geq 0}$.

\begin{Lm} \label{limitPointIsNotPeriodic}
Точка $p_{\inf}$ не является периодической точкой.
\end{Lm}
\begin{Proof}
Пусть $p_{\inf}$ - периодическая точка. Из леммы \ref{periodicComponentIsConvexPolygon}, а также, строго говоря, лемм \ref{reductionToT'}, \ref{reductionToX} следует, что существует открытый многоугольник $comp(p)$, содержащий $p$ и состоящий из периодических точек. Однако заметим, что если включить в $X$ его границу и объявить точки $\partial X$ граничными (например, зафиксировав $f$ неопределенным для этих точек), то тогда лемма \ref{typeInheritance2} будет выполнена как для $\Gamma$, так и для, очевидно, $\Gamma_1$. Следовательно, границы всех многоугольников $\Gamma_1^j(X)$, $j \in \bbz_{\geq 0}$, состоят лишь из граничных точек, а в силу того, что периметр многоугольников $\Gamma_1^j(X)$ стремится к нулю, эти границы подходят к $p_{\inf}$ бесконечно близко, что противоречит существованию $comp(p)$. Следовательно, $p_{\inf}$ не есть периодическая, QED. 
\end{Proof}

\begin{Lm} \label{limitPointIsNotBoundary}
Точка $p_{\inf}$ не является граничной.
\end{Lm}

\begin{Proof}
Пусть $p_{\inf}$ граничная. Тогда из леммы \ref{boundaryPointsAreOpenRaysAndSegments} (и, конечно, лемм \ref{reductionToT'}, \ref{reductionToZ'}, \ref{reductionToX}) следует, существует открытый отрезок $seg$, содержащий $p_{\inf}$ и состоящий лишь из граничных точек. Без ограничения общности можно считать, что длина этого отрезка равна $2l$, $l > 0$, а $p_{\inf}$ - его середина.

Рассмотрим последовательность периодических компонент $C_0$, $C_1$, $C_2$, \ldots, т.ч. $C_0 = \omega_a$, $C_1 = f(\Gamma(\omega_a))$, и $\forall k \in \bbz_{\geq 2}$: $C_k = \Gamma_1{C_{k-2}}$. Из рис. $\ref{pic:AperiodicPointSpiral-12}$ и того факта, что $p_{inf}$ лежит строго внутри $\Gamma_1(X)$ следует, что независимо от направления, как минимум один из лучей - продолжений отрезка $seg$ пересекает как минимум одну из {\it открытых} компонент $C_0$, $C_1$, причем расстояние от $p_{\inf}$ до точки пересечения не превосходит $d$, где $d = \sup\limits_{q \in C_0 \cup C_1} |p_{\inf}q|$ (здесь $|pq|$ есть длина отрезка, соединяющего точки $p,q \in \bbrr$); следовательно, $l \leq d$. С другой стороны, $\Gamma_1$ есть композиция поворота вокруг точки $p_{\inf}$ и сжатия плоскости в $\lambda = |AW_9| / |AB| < 1$ с центром в $p_{\inf}$; следовательно, один из вышеописанных лучей пересекает и одну из компонент $C_2, C_3$, причем расстояние от $p_{\inf}$ до точки пересечения не превосходит $\lambda^2 d$; отсюда получаем, что $l \leq \lambda d$. Повторяя данное рассуждение, можно показать по индукции, что $\forall k \in \mathbb{N}$: $l \leq \lambda^{2k} d$; следовательно, $l \leq 0$, что приводит нас к противоречию. 
\end{Proof}

\begin{figure}[h!]
\begin{center}
\includegraphics[width=170mm]{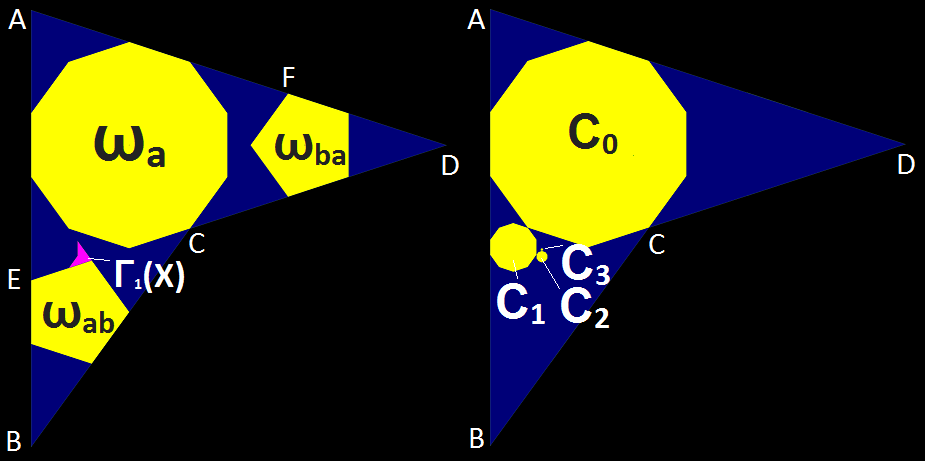}
\caption{Преобразование $\Gamma$ и начало последовательности $(C_i)$}
\label{pic:AperiodicPointSpiral-12}
\end{center}
\end{figure}

Согласно леммам \ref{limitPointIsNotPeriodic}, \ref{limitPointIsNotBoundary}, $p_{\inf}$ не может быть ни периодической, ни граничной; следовательно, $p_{\inf}$ есть апериодическая точка системы $(X, f)$, что завершает доказательство теоремы \ref{MainTreorem}.

\end{subsection}
\begin{subsection}{Доказательство теоремы 2}
Для доказательства теоремы \ref{MainTreorem2} введем пару определений.

\begin{Def}
Пусть $W \subset \bbrr$ --- многоугольник. Тогда {\it разбиением} $W$ назовем конечное множество многоугольников $\mathcal{W} = \{W_1, W_2, \ldots, W_k\}$, $k \in \bbz_+$, т.ч.:

\begin{itemize}
    \item $W = \bigcup\limits_{j=1}^{k} W_j$;
    \item $\forall i,j, 1 \leq i < j \leq k$: $\intt(W_i) \cap \intt(W_j) = \emptyset$.
\end{itemize}

\end{Def}

\begin{Def}
Пусть $W \subset \bbrr$ --- многоугольник, а $\mathcal{W}^1, \mathcal{W}^2$ - его разбиения. Тогда $\mathcal{W}^2$ является подразбиением $\mathcal{W}^1$, если $\forall Q \in \mathcal{W}^2 ~\exists P \in \mathcal{W}^1$: $Q \subset P$.
\end{Def}

Другими словами, подразбиение разбиения $\mathcal{W}^1$ есть объединение разбиений многоугольников, входящих в состав $\mathcal{W}^1$.

С разбиениями связана следующая простая лемма.

\begin{Lm} \label{figuresOfMeasure0}
Пусть $W$ есть произвольный многоугольник. Рассмотрим последовательность $(\mathcal{W}^l)_{l \in \bbz_+}$ разбиений $W$ на конечное число многоугольников, устроенную следующим образом:

\begin{enumerate}
\item $\forall l \in \bbz_+$: $\mathcal{W}^{l+1}$ есть подразбиение $\mathcal{W}^l$;

\item каждый встречающийся в хотя бы одном из разбиений последовательности $(\mathcal{W}^l)_{l \in \bbz_+}$ многоугольник в каждом разбиении покрашен в зеленый или красный цвет, причем этот цвет может зависеть от конкретного разбиения;

\item $\forall l \in \bbz_+$: если $A \in \mathcal{W}^l$ и $A$ - красный в разбиении $\mathcal{W}^l$, то $A \in \mathcal{W}^{l+1}$ и $A$ - красный в разбиении $\mathcal{W}^{l+1}$;

\item существуют числа $\epsilon \in \bbr_+$ и $k \in \bbz_+$ такие, что для любого $l \in \bbz_+$ и любого зеленого многоугольника $U \in \mathcal(W^l)$, гарантируется, что в $l+k$-ом разбиении красные фигуры, лежащие внутри $U$, обладают суммарной площадью не меньшей, чем $\epsilon A$, где $A$ есть площадь фигуры $U$.
\end{enumerate}

Тогда объединение участвующих в разбиениях красных фигур образует в $W$ множество полной меры.
\end{Lm}

\begin{Proof}
Обозначим за $Area(U)$ площадь произвольного многоугольника $U$. 
Пусть $RedArea(l)$, $l \in \bbz_+$, есть суммарная площадь красных фигур, участвующих в разбиении $\mathcal{W}^l$; аналогично определим $GreenArea(l)$ (очевидно, что $\forall l \in \bbz_+: RedArea(l) + GreenArea(l) = Area(W)$). Из условия 3 следует, что последовательность $(RedArea(l))$ неубывает с ростом $l$, а $(GreenArea(l))$ невозрастает. С другой стороны, из условия 4 следует, что $\forall l \in \bbz_+$: $GreenArea(l + k) \leq (1 - \epsilon) * GreenArea(l)$. Так как $\epsilon > 0$, то $GreenArea(l) \underset{l \rightarrow +\infty}{\rightarrow} 0$. Это означает, что $\forall \delta \in \bbr_+$: все 
точки многоугольника $W$, которые не лежат ни в одном из красных многоугольников, могут быть покрыты зелеными многоугольниками разбиения $\mathcal{W}^{l_{\delta}}$, где $l_{\delta}$ таково, что суммарная площадь этих многоугольников не превосходит $\delta$. Таким образом, дополнение в $W$ к объединению всех красных фигур образует множество меры нуль, QED.
\end{Proof}

С помощью леммы \ref{figuresOfMeasure0} докажем полноту меры периодических точек в системе $(X, f)$; этого, как следует из лемм \ref{reductionToX}, \ref{reductionToZ'}, \ref{reductionToT'} будет достаточно для доказательства теоремы \ref{MainTreorem2}.

\begin{Lm} \label{periodicPointsOfXAreOfFullMeasure}
Периодические относительно преобразования $f$ точки образуют в $X$ множество полной меры.
\end{Lm}

\begin{Proof}
Используем лемму \ref{figuresOfMeasure0}. Рассмотрим последовательность раскрашенных в красный и зеленый цвета разбиений $(\mathcal{X}^l)_{l \in \bbz_+}$ многоугольника $X$, устроенную следующим образом:

\begin{itemize}
    \item $(\mathcal{X}^1)$ состоит из двух красных треугольников $ABE$ и $CED$;
    \item $\forall l \in \bbz_+$: $\mathcal{X}^{l+1} = \{\bigcup\{\overline{f^j(\Gamma(\intt(P)))}\} | P \in \mathcal{X}^l, 0 \leq j < t(P)\} \cup \{\overline{\omega_a}, \overline{\omega_{ba}}, \overline{\omega_{ab}} \}$; здесь $t(P)$ есть минимальное такое $t \in \bbz_+$, что $f_t(\intt(P)) \subset \Gamma(X)$;
    \item в терминах предыдущего пункта, многоугольник $\overline{f^j(\Gamma((\intt(P)))}$ зеленый в разбиении $(\mathcal{X})^{l+1}$, если и только если $P$ зеленое в разбиении $(\mathcal{X})^l$;
    \item фигуры $\overline{\omega_a}, \overline{\omega_{ba}}, \overline{\omega_{ab}}$ являются красными во всех разбиениях, в которых участвуют. 
\end{itemize}

Из рис. \ref{pic:RedGreenPETs} очевидно следует, что $\mathcal{X}^2$ есть разбиение $X$; тот же факт, что $\mathcal{X}^l$ является разбиением для любого $l \in \bbz$, легко доказать по индукции. Более того, на основе леммы \ref{selfSimilarityX} также можно вывести по индукции, что в разбиении $\mathcal{X}^l$, $l \in \bbz_+$, все красные многоугольники являются периодическими компонентами (точнее, их замыканиями), а зеленые являются траекториями первого возвращения в $\Gamma^{l-1}(X)$ относительно $f$ открытых треугольников $\Gamma_{l-1}(\intt(ABE))$ и $\Gamma_{l-1}(\intt(CED))$.

\begin{figure}[h!]
\begin{center}
\includegraphics[width=170mm]{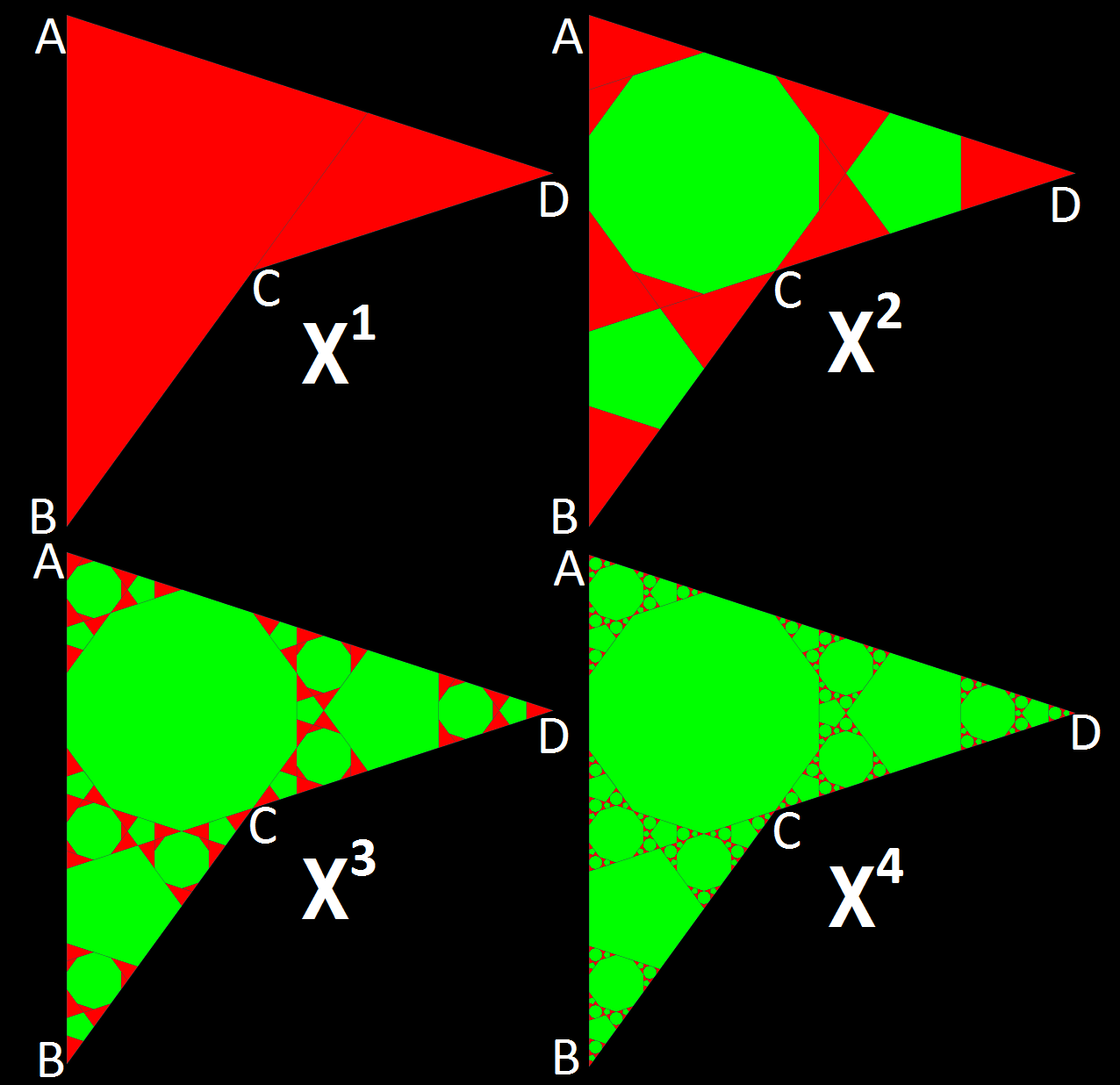}
\caption{Разбиения $\mathcal{X}^1$, $\mathcal{X}^2$, $\mathcal{X}^3$, $\mathcal{X}^4$}
\label{pic:RedGreenPETs}
\end{center}
\end{figure}

Из вышесказанного очевидно, что выполнены условия 1-3 леммы \ref{figuresOfMeasure0}. Условие 4 той же леммы также выполнено; это следует из того, что при переходе от $\mathcal{X}^l$ к $\mathcal{X}^{l+1}$, $l \in \bbz_+$, верно:

\begin{itemize}
    \item $\forall P \in \mathcal{X}^l$ т.ч. $P$ - зеленый в $\mathcal{X}^l$, и $P = \overline{f^j(\Gamma^{l-1}(\intt(ABE))}$, $j \in \bbz_{\geq 0}$: многоугольники $P_a = \overline{f^j(\Gamma^{l}(\omega_a))}$ и $P_{ab} = \overline{f^j(\Gamma^{l}(\omega_{ab}))}$ таковы, что:
    \begin{itemize}
        \item $P_a, P_{ab} \subset P$;
        \item $P_a, P_{ab}$ - красные многоугольники в $\mathcal{X}^{l+1}$;
        \item $\frac{Area(P_a \cup P_{ab})}{Area(P)} = \frac{Area(\omega_{a} \cup \omega_{ab})}{Area(ABE)}$;
    \end{itemize}
    \item $\forall P \in \mathcal{X}^l$ т.ч. $P$ - зеленый в $\mathcal{X}^l$, и $P = \overline{f^j(\Gamma^{l-1}(\intt(CED))}$, $j \in \bbz_{\geq 0}$: многоугольник $P_{ba} = \overline{f^j(\Gamma^{l}(\omega_{ba}))}$ таков, что:
    \begin{itemize}
        \item $P_{ba} \subset P$;
        \item $P_{ba}$ - красные многоугольники в $\mathcal{X}^{l+1}$;
        \item $\frac{Area(P_{ba})}{Area(P)} = \frac{Area(\omega_{ba})}{Area(CED)}$.
    \end{itemize}
    
\end{itemize}

Таким образом, условие 4 леммы \ref{figuresOfMeasure0} выполнено, если взять $k = 1$ и $\epsilon = \min(\frac{Area(\omega_{a} \cup \omega_{ab})}{Area(ABE)}, \frac{Area(\omega_{ba})}{Area(CED)})$. Следовательно, все условия леммы \ref{figuresOfMeasure0}, а согласно этой лемме, красные многоугольники, являющиеся замыканиями периодических компонент, образуют в $X$ множество полной меры. Строго говоря, периодические компоненты - это открытые многоугольники, отличающиеся от своих замыканий; однако тот факт, что объединение замыканий счетного числа открытых многоугольников отличается от объединения счетного числа тех же открытых многоугольников лишь множеством точек меры нуль (ибо это счетное число отрезков), завершает доказательство леммы. 
\end{Proof}

\end{subsection}
\begin{subsection}{Нахождение абелизаций кодов периодов периодических компонент для $(X, f)$}.
В следующих нескольких разделах мы найдем множество периодов точек для внешнего биллиарда $T$ вне правильного десятиугольника $\gamma$. Мы будем действовать по следующему плану. В этом разделе, мы найдем множество 
$C_X = \{c(\rho_{Xper}(p)) | p \subset X, p \text{ - периодическая компонента}\}$, где $c(w): \{a,b\}^* \rightarrow \bbz_{\geq 0}^2$ есть гомоморфизм абелизации, определенный аналогично определению \ref{def:abelization}. После этого, с помощью леммы \ref{reductionToX} мы вычислим аналогичным образом определенное множество $C_{Z' \cup \beta_4}$. На следующем шаге, по лемме \ref{fromZ'PeriodCodesToV'PeriodCodesAbel} мы найдем $C_{V'}$; затем, по лемме \ref{fromIPStructureToPeriodAbel} вычислим и множество всевозможных периодов компонент $B$ для $T$ вне $\Gamma$. Наконец, согласно лемме \ref{fixedPeriods}, преобразуем множество $B$ в $B_2$, просто добавив в $B$ удвоения всех нечетных чисел в $B$.

Пусть $\sigma: \{a, b\} \rightarrow \{a, b\}^*$ - подстановка, т.ч. $\sigma(a) = aababaa$, $\sigma(b) = aaa$. Из рис. \ref{pic:Xf-SelfSimilarity} очевидно доказательство следующей леммы, являющейся уточнением третьего пункта леммы \ref{selfSimilarityX}.

\begin{Lm} \label{selfSimilarityXSubstitution}
Пусть $p \subset X$ - периодическая компонента, и пусть $q = \Gamma(p)$. Тогда $\rho_{Xper}(q) = \sigma(\rho_{Xper}(p))$.
\end{Lm}

Очевидно, лемма \ref{selfSimilarityXSubstitution} могла бы быть сформулирована и для периодической или апериодической точки $p \in X$ и ее кода $\rho_X(p)$; однако для наших целей достаточно ограничиться лишь периодическими компонентами.

Из леммы \ref{genPeriodicStrong2} и того факта, что если $p \subset X$ - периодическая компонента, то и $f(p)$ - периодическая компонента, причем $c(\rho_{Xper}(p)) = c(\rho_{Xper}(f(p)))$ следует, что $C_X = \{c(\rho_{Xper}(\Gamma^l(\omega))) | l \in \bbz_{\geq 0}, \omega \in \{\omega_a, \omega_{ba}\}\}$.

Пусть $M_{\sigma}$ - матрица подстановки $\sigma$. Согласно определению $\sigma$, $M_{\sigma} = \begin{pmatrix}5 & 3 \\ 2 & 0 \end{pmatrix}$, и если $p \subset X$ - периодическая компонента, то $c(\rho_{Xper}(\Gamma(p))) = M_{\sigma}c(\rho_{Xper}(p))$. Следовательно, $C_X = \{M_{\sigma}^lc(\rho_{Xper}(\omega)) | \omega \in \{\omega_a, \omega_{ba}\}\} = \{M_{\sigma}^lg | g \in \{ \vectwo{1}{0}, \vectwo{1}{1}\}\}$ (ибо $c(\omega_a) = \vectwo{1}{0}$, $c(\omega_{ba}) = \vectwo{1}{1}$ ).

С другой стороны, легко получить с помощью линейной алгебры и доказать по индукции, что $\forall x, y \in \bbr, l \in \bbz_{\geq 0}$: $M_{\sigma}^l\vectwo{x}{y} = \frac{x-3y}{7} (-1)^l \vectwo{1}{-2} + \frac{y+2x}{7} 6^l \vectwo{3}{1}$. Подставляя в качестве $\vectwo{x}{y}$ векторы $\vectwo{1}{0}, \vectwo{1}{1}$, получаем следующую лемму.

\begin{Lm} \label{abelizationsOfX}
$C_X = \{\frac{1}{7}\pvectwo{6*6^l + (-1)^l}{2*6^l - 2*(-1)^l}, \frac{1}{7}\pvectwo{9*6^l - 2*(-1)^l}{3*6^l + 4*(-1)^l} | l \in \bbz_{\geq 0}\}$.
\end{Lm}

\end{subsection}
\begin{subsection}{Нахождение абелизаций кодов периодических компонент для $(V', T')$}.
Вычислим $C_{Z' \cup \beta_4}$ и $C_{V'}$. С помощью леммы \ref{reductionToX}, напрямую получаем следующую лемму.

\begin{Lm} \label{abelizationOfZ'AndBeta4}
$C_{Z' \cup \beta_4} = \{\frac{1}{7}\vecfive{12*6^l+2*(-1)^l}{24*6^l-3*(-1)^l}{0}{0}{0}, 
                         \frac{1}{7}\vecfive{18*6^l-4*(-1)^l}{36*6^l+6*(-1)^l}{0}{0}{0}, \\
                         \frac{1}{7}\vecfive{0}{0}{6*6^l+(-1)^l}{8*6^l-(-1)^l}{0},
                         \frac{1}{7}\vecfive{0}{0}{9*6^l-2*(-1)^l}{12*6^l+2*(-1)^l}{0},
                         \vecfive{1}{0}{0}{0}{0}, \vecfive{0}{1}{0}{0}{0},
                         \vecfive{0}{0}{1}{0}{0}, \vecfive{0}{0}{0}{1}{0} | l \in \bbz_{\geq 0}
                       \}$.
\end{Lm}

Чтобы вычислить $C_{V'}$, воспользуемся леммой \ref{fromZ'PeriodCodesToV'PeriodCodesAbel}, согласно которой, каждый вектор из $C_{Z' \cup \beta_4}$ нужно домножить слева на матрицу $M_{\psi}$ произвольное число раз. Пусть $C^1_{V'} = \{M_{\psi}w~|~w \in C_{Z' \cup \beta_4} \}$. Из устройства $M_{\psi}$ следует, что первые три компоненты векторов $C^1_{V'}$ нулевые. В силу этого, введем преобразования $pr_{25}: \bbr^2 \rightarrow \bbr^5$, т.ч. $\forall u,v \in \bbr$: $pr_{25}\vectwo{u}{v} = \vecfive{0}{0}{0}{u}{v}$, и $pr_{52} = pr_{25}^{-1}$; пусть $C^1_{V'2} = pr_{52}(C^1_{V'})$.

Также введем вспомогательную матрицу $M_{\psi1} = \begin{smallmatrix}1 & 0 \\ 1 & 1 \end{smallmatrix}$. Очевидно, коэффициенты подобраны таким образом, что $\forall w \in \bbrr$: $M_{\psi}pr_{25}(w) = pr_{25}(M_{\psi1}w)$.

Используя вышесказанное, лемму \ref{fromZ'PeriodCodesToV'PeriodCodesAbel} можно записать в виде:

$C_{V'} = C_{Z' \cup \beta_4} \cup \{pr_{25}(M_{\psi1}^kw~)|~w \in  C^1_{V'2}, k \in \bbz_{\geq 0} \}$.

Из определения легко получить, что

$C^1_{V'2} = \{\frac{1}{7}\vectwo{120*6^l-(-1)^l}{36*6^l-(-1)^l},
               \frac{1}{7}\vectwo{180*6^l+2*(-1)^l}{54*6^l+2*(-1)^l},
                \frac{1}{7}\vectwo{20*6^l+(-1)^l}{14*6^l},
                \frac{1}{7}\vectwo{30*6^l-2*(-1)^l}{21*6^l}, \\
                \vectwo{4}{1}, \vectwo{3}{1}, \vectwo{2}{1}, \vectwo{1}{1} 
             \}$.

С другой стороны, легко доказать по индукции по $k$, что $\forall x,y \in \bbr, k \in \bbz_{\geq 0}$: $M^k_{\psi1}\vectwo{x}{y} = \vectwo{x}{y+kx}$. Отсюда следует

\begin{Lm} \label{abelizationsOfV'}
$C_{V'} = C_{V'f} + pr_{25}(C_{V'o})$, где $C_{V'f} = C_{Z' \cup \beta_4}$,

$C_{V'o} = \{\frac{1}{7}\vectwo{120*6^l-(-1)^l}{(36+120k)*6^l-(k+1)*(-1)^l},
             \frac{1}{7}\vectwo{180*6^l+2*(-1)^l}{(54+180k)*6^l+2*(k+1)*(-1)^l}, \\
             \frac{1}{7}\vectwo{20*6^l+(-1)^l}{(14+20k)*6^l+k*(-1)^l},
             \frac{1}{7}\vectwo{30*6^l-2*(-1)^l}{(21+30k)*6^l-2k*(-1)^l},
             \vectwo{4}{1+4k}, \vectwo{3}{1+3k}, \vectwo{2}{1+2k}, \\ \vectwo{1}{1+k}
             |k, l \in \bbz_{\geq 0}\}$.
\end{Lm}

\end{subsection}
\begin{subsection} {Нахождение периодов и доказательство теоремы \ref{MainTheorem3}}
Чтобы найти множество $B_c$ периодов периодических компонент внешнего биллиарда вне правильного десятиугольника, достаточно применить лемму \ref{fromIPStructureToPeriodAbel} к каждому из столбцов $w \in C_{V'}$ и добавить в $B_c$, изначально пустое, число $r(w) = \frac{10s(w)}{\rgcd(10, t(w))}$, где $s(w) = (1,1,1,1,1)*w$, $t(w) = (1,2,3,4,5)*w$. Найдем $s(w)$ и $t(w)$ для всех серий столбцов, встречающихся в $C_{V'}$, согласно леммам \ref{abelizationOfZ'AndBeta4}, \ref{abelizationOfV'}. При подсчете мы используем следующие простые факты теории чисел:

\begin{itemize}
    \item $\forall a, b \in \bbz_+$: $\rgcd(a, b) = \rgcd(7a, b)$; следовательно, $\forall w \in \bbr^5$: $\rgcd(10, t(w)) = \rgcd(10, t(7w))$. Этот факт позволит нам игнорировать множитель $\frac{1}{7}$ при подсчете НОД-а;

    \item $6^0 \equiv_{10} 1$, и $\forall l \in \bbz_+$: $6^l \equiv_7 6$.
    
\end{itemize}
Итак,
\begin{itemize}
    \item пусть $w = \frac{1}{7}\vecfive{12*6^l+2*(-1)^l}{24*6^l-3*(-1)^l}{0}{0}{0}$, $l \in \bbz_{\geq 0}$. Тогда $7s(w) = 6^{l+2}-(-1)^l$, $7t(w) = (10*6^{l+1} - 4*(-1)^l)$, $\rgcd(10, t(w)) = \rgcd(10,  - 4(-1)^l) = \rgcd(10, 4) = 2$. Следовательно, $r(w) = \frac{5}{7}(6^{l+2}-(-1)^l)$. 

    \item Пусть $w = \frac{1}{7}\vecfive{18*6^l-4*(-1)^l}{36*6^l+6*(-1)^l}{0}{0}{0}$, $l \in \bbz_{\geq 0}$. Тогда $7s(w) = 9*6^{l+1} + 2*(-1)^l$, $7t(w) = (15*6^{l+1} + 8*(-1)^l)$, $\rgcd(10, t(w)) = \rgcd(10, 90 + 8(-1)^l) = \rgcd(10, 8(-1)^l) = \rgcd(10, 8) = 2$; следовательно, $r(w) = \frac{5}{7}(9*6^{l+1} + 2*(-1)^l)$.
    
    \item Пусть $w = \frac{1}{7}\vecfive{0}{0}{6*6^l+(-1)^l}{8*6^l-(-1)^l}{0}$, $l \in \bbz_{\geq 0}$. Тогда $7s(w) = 14*6^l$, $7t(w) = 50*6^l - (-1)^l$, $\rgcd(10, t(w)) = \rgcd(10, -(-1)^l) = 1$; следовательно, $r(w) = 20*6^l$.
    
    \item Пусть $w = \frac{1}{7}\vecfive{0}{0}{9*6^l-2*(-1)^l}{12*6^l+2*(-1)^l}{0}$, $l \in \bbz_{\geq 0}$. Тогда $7s(w) = 21*6^l$, $7t(w) = 75*6^l + 2*(-1)^l$. Рассмотрим два случая:
    \begin{itemize}
        \item $l = 0$; тогда $7t(w) = 77$, $\rgcd(10, t(w)) = 1$, и $r(w) = 30$;
        \item $l > 0$; тогда $\rgcd(10, t(w)) = \rgcd(10, 30*6^{l-1} + 2*(-1)^l) = \rgcd(10, 2*(-1)^l) = 2$, и $r(w) = 15*6^l$.
    \end{itemize}
    
    \item Пусть $w = \vecfive{1}{0}{0}{0}{0}$; тогда $s(w) = 1$, $t(w) = 1$, и $r(w) = 10$.
    \item Пусть $w = \vecfive{0}{1}{0}{0}{0}$; тогда $s(w) = 1$, $t(w) = 2$, и $r(w) = 5$.
    \item Пусть $w = \vecfive{0}{0}{1}{0}{0}$; тогда $s(w) = 1$, $t(w) = 3$, и $r(w) = 10$.
    \item Пусть $w = \vecfive{0}{0}{0}{1}{0}$; тогда $s(w) = 1$, $t(w) = 4$, и $r(w) = 5$.

    \item Пусть $w = pr_{25}(\frac{1}{7}\vectwo{120*6^l-(-1)^l}{(36+120k)*6^l-(k+1)*(-1)^l})$, $k,l \in \bbz_{\geq 0}$. Тогда $7s(w) = (156+120k)*6^l - (k+2)*(-1)^l$, $7t(w) = (660+600k)*6^l - (5k+9)*(-1)^l$, $\rgcd(10, t(w)) = (10, 5k+9)$. Рассмотрим два случая:
    \begin{itemize}
        \item $k = 2m, m \in \bbz_{\geq 0}$; тогда $\rgcd(10, t(w)) = 1$, и $r(w) = \frac{10}{7}((156+120*2m)*6^l - (2m+2)*(-1)^l) = \frac{20}{7}((78+120m)*6^l - (m+1)*(-1)^l)$;
        \item $k = 2m+1, m \in \bbz_{\geq 0}$; тогда $\rgcd(10, t(w)) = 2$, и $r(w) = \frac{5}{7}((276+240m)*6^l - (2m+3)*(-1)^l)$. 
    \end{itemize}

    \item Пусть $w = pr_{25}(\frac{1}{7}\vectwo{180*6^l+2*(-1)^l}{(54+180k)*6^l+2*(k+1)*(-1)^l})$,  $k,l \in \bbz_{\geq 0}$. Тогда $7s(w) = (234+180k)*6^l + (2k+4)*(-1)^l$, $7t(w) = (990+900k)*6^l + (10k+18)*(-1)^l$, $\rgcd(10, t(w)) = \rgcd(10, 18) = 2$; следовательно, $r(w) = \frac{5}{7}((234+180k)*6^l + (2k+4)*(-1)^l)$.
    
    \item Пусть $w = pr_{25}(\frac{1}{7}\vectwo{20*6^l+(-1)^l}{(14+20k)*6^l+k*(-1)^l})$. $k,l \in \bbz_{\geq 0}$. Тогда $7s(w) = (34 + 20k)*6^l + (k+1)*(-1)^l$, $7t(w) = (150+100k)*6^l+(5k+4)*(-1)^l$, $\rgcd(10, t(w)) = \rgcd(10, 5k+4)$. Рассмотрим два случая:
    \begin{itemize}
        \item $k = 2m, m \in \bbz_{\geq 0}$; тогда $\rgcd(10, t(w)) = 2$, и $r(w) = \frac{5}{7}((34 + 40m)*6^l + (2m+1)*(-1)^l)$;
        \item $k = 2m+1, m \in \bbz_{\geq 0}$; тогда $\rgcd(10, t(w)) = 1$, и $r(w) = \frac{10}{7}((20 + 40m)*6^l + (2m+2)*(-1)^l)$.
    \end{itemize}

    \item Пусть $w = pr_{25}(\frac{1}{7}\vectwo{30*6^l-2*(-1)^l}{(21+30k)*6^l-2k*(-1)^l})$. $k,l \in \bbz_{\geq 0}$. Тогда $7s(w) = (51 + 30k)*6^l - (2k+2)*(-1)^l$, $7t(w) = (225+150k)*6^l-(10k+8)*(-1)^l$, $\rgcd(10, t(w)) = \rgcd(10, 5*6^l + 8*(-1)^l)$.
    Рассмотрим два случая:
    \begin{itemize}
        \item $l = 0$; тогда $\rgcd(10, t(w)) = 1$, и $r(w) = \frac{10}{7}(28k+49) = 40k+70$;
        \item $l > 0$; тогда $\rgcd(10, t(w)) = 2$, и $r(w) = \frac{5}{7}((51 + 30k)*6^l - (2k+2)*(-1)^l)$.
    \end{itemize}
% это была последняя сложная серия

    \item Пусть $w = pr_{25}(\vectwo{4}{1+4k})$, $k \geq 0$; тогда $s(w) = 5+4k$, $t(w) = 20k+21$, $\rgcd(10, t(w)) = 1$; следовательно, $r(w) = 40k+50$.
    
    \item Пусть $w = pr_{25}(\vectwo{3}{1+3k})$, $k \geq 0$; тогда $s(w) = 4+3k$, $t(w) = 15k+17$, $\rgcd(10, t(w)) = \rgcd(10, 5k+7)$. Рассмотрим два случая:
    \begin{itemize}
        \item $k = 2m$, $m \geq 0$; тогда $\rgcd(10, t(w)) = 1$, и $r(w) = 30k+40 = 60m + 40$;
        \item $k = 2m+1$, $m \geq 0$; тогда $\rgcd(10, t(w)) = 2$, и $r(w) = 15k+20 = 30m + 35$.
    \end{itemize}
    
    \item Пусть $w = pr_{25}(\vectwo{2}{1+2k})$, $k \geq 0$; тогда $s(w) = 3+2k$, $t(w) = 10k+13$, $\rgcd(10, t(w)) = 1$; следовательно, $r(w) = 20k+30$.
    
    \item Пусть $w = pr_{25}(\vectwo{1}{1+k})$, $k \geq 0$; тогда $s(w) = 2+k$, $t(w) = 5k+9$, $\rgcd(10, t(w)) = \rgcd(10, 5k+9)$. Рассмотрим два случая:
    \begin{itemize}
        \item $k = 2m$, $m \geq 0$; тогда $\rgcd(10, t(w)) = 1$, и $r(w) = 10k+20 = 20m + 20$;
        \item $k = 2m+1$, $m \geq 0$; тогда $\rgcd(10, t(w)) = 2$, и $r(w) = 5k+10 = 10m + 15$.
    \end{itemize}
    
\end{itemize}

Объединив полученные числа и серии $r(w)$ и заметив, что числа каждой из серий либо все четные, либо все нечетные, получаем следующую лемму, являющуюся уточненной версией теоремы \ref{MainTheorem3}.

\begin{Lm}
Пусть $B = \{\frac{5}{7}(6^{l+2}-(-1)^l), \frac{5}{7}(9*6^{l+1} + 2*(-1)^l), 20*6^l,
          30, 90*6^l, 10, 5, \frac{20}{7}((78+120k)*6^l - (k+1)*(-1)^l),
          \frac{5}{7}((276+240k)*6^l - (2m+3)*(-1)^l), \frac{5}{7}((234+180k)*6^l + (2k+4)*(-1)^l),
          \frac{5}{7}((34 + 40k)*6^l + (2k+1)*(-1)^l), \frac{10}{7}((20 + 40k)*6^l + (2k+2)*(-1)^l),
          40k+70, \frac{5}{7}((306 + 180k)*6^l + (2k+2)*(-1)^l), 40k+50, 60k + 40, 30k + 35,
          20k+30, 20k + 20, 10k + 15|k, l \in \bbz_{\geq 0}\}$,
         
      $B_2 = B \cup \{\frac{10}{7}(6^{l+2}-(-1)^l),  \frac{10}{7}((276+240k)*6^l - (2k+3)*(-1)^l), \frac{10}{7}((34 + 40k)*6^l + (2k+1)*(-1)^l), 60k+70, 20k+30| k, l \in \bbz_{\geq 0}\}$.
      
Тогда $B$ и $B_2$ суть множества периодов периодических компонент и периодических точек для внешнего биллиарда вне правильного десятиугольника.
\end{Lm}

Утверждение леммы про $B_2$ получено с помощью леммы \ref{fixedPeriods}.

\end{subsection}

\end{section}

\begin{section}{Заключение}
В данной работе было проведено полное исследование внешнего биллиарда вне правильного десятиугольника, завершающее <<программу Шварца>>. Основными наблюдениями, на которых базируются доказательства теорем, являются леммы \ref{selfSimilarityBig} и \ref{selfSimilarityX}, устанавливающие наличие самоподобных структур в множествах периодических/апериодических точек вне стола $\gamma$. Самоподобие, возникающее в первой из этих лемм, не случайно; аналогичное свойство можно доказать для внешнего биллиарда вне правильного $n$-угольника для произвольного четного $n$; в случае же нечетного $n$, можно установить соответствие между периодическими/апериодическими структурами вне правильных $n$- и $2n$-угольника (для случая $n = 5$, такое соответствие описано в \cite{BC11}). Однако самоподобие, задаваемое леммой \ref{selfSimilarityX}, таким общим свойством не является; похожие самоподобия обнаружены лишь для случаев $n = 5, 10, 8, 12$, причем в последнем случае были задействованы доказательные компьютерные вычисления. Более того, проведенные Р.Шварцем компьютерные эксперименты \cite{IJJ15} показали, что периодические структуры для внешнего биллиарда вне правильного семиугольника обладают существенно более сложной структурой, нежели случаи $n = 5,10,8,12$; в частности, в этих случаях периодическая компонента может быть неравносторонним многоугольником. Так или иначе, проблемы периодичности для правильных $n$-угольников, где $n \notin \{3,4,6,5,10,8,12\}$, остаются открытыми.

\end{section}

%\bibliographystyle{utf8gost705u}
%\bibliography{biblio}
%\end{thebibliography}

\end{document}